\documentclass[a4paper,12pt,twoside,draft]{article}
\usepackage{amsmath,amssymb,amsthm}
\usepackage{upref}

\evensidemargin = \oddsidemargin

\usepackage[hyperindex,bookmarksnumbered,
           ]{hyperref}
\theoremstyle{plain}

  \newtheorem{theorem}{Theorem}[section]
  \newtheorem{prop}{Proposition}[section]
  \newtheorem{lemma}{Lemma}[section]
  \newtheorem{cor}{Corollary}[section]
  \newtheorem{remark}{Remark}[section]
  \newtheorem{definition}{Definition}[section]

\numberwithin{equation}{section}

\newcommand{\bo}{{\hfill\loota}}
\newcommand{\loota}{\hbox{\enspace{\vrule height 7pt depth 0pt width
      7pt}}}
\newcommand{\bp}{\mbox{\boldmath $\psi$}}

\newcommand{\R}{\mathbb{R}}

\newcommand{\C}{\mathbb{C}}

\newcommand{\al}{\alpha}
\newcommand{\ov}{\overline}

\newcommand{\ve}{\varepsilon}

\newcommand{\om}{\omega}
\newcommand{\Om}{\Omega}

\newcommand{\na}{\nabla}
\newcommand{\si}{\sigma}
\newcommand{\Si}{\Sigma}

\newcommand{\ds}{\displaystyle}
\newcommand{\la}{\label}

\newcommand{\beqn}{\begin{eqnarray}}
\newcommand{\eeqn}{\end{eqnarray}}
\newcommand{\be}{\begin{equation}}
\newcommand{\ee}{\end{equation}}

\newcommand{\pa}{\partial}
\newcommand{\fr}{\frac}

\newcommand{\ci}{\cite}
\newcommand{\re}{\ref}

\renewcommand{\Re}{\mathop{\mathrm{Re}}} 
\renewcommand{\Im}{\mathop{\mathrm{Im}}} 

\DeclareMathOperator{\const}{const}
\DeclareMathOperator{\supp}{supp}

\title{\footnotesize{\textit{}}\hfill\hbox{}
\newline\hbox{}\newline\Large{}}

\begin{document}

\begin{titlepage}

\hspace{12cm}

\begin{center}
{\Large\bf
Scattering of Solitons
for Schr\"odinger \medskip\\
Equation
Coupled to a Particle
}\\
\vspace{1cm}
{\large Alexander Komech}
\footnote{
On leave Department Mechanics and Mathematics
of Moscow State University.
Supported partly by
Max-Planck Institute of Mathematics in the Sciences (Leipzig),
 and Wolfgang Pauli Institute of Vienna University.}
\medskip\\
{\it Faculty of Mathematics of
Vienna University\\
Nordbergstrasse 15,
1090 Vienna, Austria
}\\
email: komech@mathematik.tu-muenchen.de
\bigskip\\
 {\large Elena Kopylova}\footnote{
}
\medskip\\
{\it
M.V.Keldysh Institute of Applied Mathematics RAS\\
Miusskaya sq. 4,
 125047 Moscow, Russia
}\\
email: ek@vpti.vladimir.ru

\end{center}

\vspace{2cm}
\begin{abstract}
We establish soliton-like asymptotics for finite energy solutions to
the Schr\"odinger equation coupled to a nonrelativistic classical particle.
Any solution with initial state close to the solitary
manifold, converges to a sum of traveling wave and outgoing free wave.
The convergence holds in global energy norm. The proof uses
spectral theory and the symplectic projection onto solitary
manifold in the Hilbert phase space.

\end{abstract}
\end{titlepage}


\setcounter{equation}{0}

\section{Introduction}

We continue the study of coupled systems of wave fields and
particles. In \cite{ikv}  the Klein-Gordon equation
coupled to a relativistic particle has been considered. Here we extend the result to the
Schr\"odinger equation coupled to a nonrelativistic particle. We prove
the long time convergence to the sum of a soliton and dispersive wave.
The convergence holds in global energy norm for finite energy solution
with initial state close to the solitary manifold.

We consider the  Schr\"odinger wave function $\psi(x)$ in $\R^3$, coupled to a nonrelativistic
particle with position $q$ and momentum $p$, governed by
\begin{equation} \label{S}
 \left\{\begin{array}{l}
      i\dot\psi(x,t)=-\Delta\psi(x,t)+m^2\psi(x,t)  + \rho(x-q(t))
  \\
    \ddot q(t)=\ds\frac 12\ds\int\bigl[\overline{\psi}(x,t)\nabla\rho(x-q(t))
    +\psi(x,t)\nabla\ov\rho(x-q(t))\bigr]dx,
\end{array}\right|\; x\in\R^3
\end{equation}
where $m>0$.
Denote  $\psi_1=\Re\psi,\;\psi_2=\Im\psi,\;\rho_1=\Re\rho,\;\rho_2=\Im\rho$.\\
Then the system \eqref{S} becomes
\begin{equation} \label{SS}
 \left\{\begin{array}{l}
      \dot\psi_1(x,t)=-\Delta\psi_2(x,t)+m^2\psi_2(x,t) + \rho_2(x-q(t)),\\
      \dot\psi_2(x,t)=\Delta\psi_1(x,t)-m^2\psi_1(x,t) - \rho_1(x-q(t)),\;
  \\
    \ddot q(t)=\ds\int\bigl(\psi_1(x,t)\nabla\rho_1(x-q(t))
    +\psi_2(x,t)\nabla\rho_2(x-q(t))\bigr)dx.
\end{array}\right|\; x\in\R^3
\end{equation}
This is a Hamilton system with the Hamilton functional
\begin{multline}\label{Hamil}
  {\cal H}(\psi_1,\psi_2,q,\dot q)=\frac 12\int\Big(
  |\nabla \psi_1(x)|^2+|\nabla \psi_2(x)|^2+m^2|\psi_1(x)|^2+m^2|\psi_2(x)|^2\Big)~dx\\
  +\int\Big(\psi_1(x)\rho_1 (x-q)+\psi_2(x)\rho_2 (x-q)\Big)dx+\frac 12|\dot q|^2.
\end{multline}
We consider the Cauchy problem for the Hamilton system (\re{SS}) which we write as
\be\la{WP2.1}
  \dot Y(t)=F(Y(t)),\quad t\in\R;\quad Y(0)=Y_0.
\ee
Here $Y(t)=(\psi_1(t), \psi_2(t), q(t), p(t))$, $p(t):=\dot q(t)$,
$Y_0=(\psi_{01},\psi_{02}, q_0, p_0)$, and all
derivatives are understood in the sense of distributions.
Below we always deal with column vectors but often write them as row vectors.
The system \eqref{SS} is translation-invariant and admits soliton solutions
\be\la{sosol}
  Y_{a,v}(t)=(\psi_{v1}(x-vt-a),\psi_{v2}(x-vt-a), vt+a,v),
\ee
for all $a,v\in\R^3$ with $|v|<2m$.
The states $S_{a,v}:=Y_{a,v}(0)$ form the solitary manifold
\be\la{soman}
{\cal S}:=\{ S_{a,v}: a,v\in\R^3, |v|<2m \}.
\ee
Our main result is the soliton asymptotics of type
\be\la{Solas}
  \psi(x,t)\sim\psi_{v_\pm }(x-v_\pm t-a_\pm)+W_0(t){\bp}_\pm,\quad t\to\pm\infty,
\ee
for solutions to  (\re{S}) with initial data close to the solitary manifold ${\cal S}$.
Here $\psi_{v_\pm }=\psi_{v_\pm 1}+i\psi_{v_\pm 2}$, $W_0(t)$ is the dynamical group of
the free Schr\"odinger equation, $\bp_\pm$ are the corresponding
{\it asymptotic scattering states}, and the asymptotics hold
{\it in the global energy norm}, i.e. in the norm of the Sobolev space $H^1(\R^3)$.
For the particle trajectory we prove that
\be\la{qq}
  \dot q(t)\to v_\pm,\quad q(t) \sim v_\pm t+a_\pm,\quad t\to\pm\infty.
\ee
The results are established under the following conditions
on the complex valued charge distributions $\rho$:
\be\la{ro}
(1+|x|)^\beta\rho,\quad (1+|x|)^\beta\nabla\rho,\quad (1+|x|)^\beta\nabla\nabla\rho\in L^2(\R^3),
\ee
with some $\beta>3/2$.
We require that all ``modes'' of the wave field are coupled to the particle,
this is formalized by the Wiener condition
\be\la{W}
  \hat\rho(k)=(2\pi)^{-3/2}\int\limits \, e^{ik x}\rho(x)dx\not=0
  \mbox{ \,\,\,for\,\,all\,\, }k\in\R^3\,.
\ee
It is an analogue of the Fermi Golden Rule: the coupling term  $\rho(x-q)$ is not orthogonal
to the eigenfunctions $e^{ikx}$ of the continuous spectrum of the linear part
of the   equation (cf. \ci{BP3,Sig,SW3,SW4}).

Similar results were proved for the first time by Buslaev and Perelman
\cite{BP1,BP2} for 1D translation invariant Schr\"odinger equation,
and extended by Cuccagna \cite{Cu} for nD case, $n\ge 3$.
In \cite{ikv}  the Klein-Gordon equation coupled to
a particle, is considered.

For the proofs of the asymptotics
(\re{Solas}) and (\re{qq}), we develop the approach \cite{ikv} based
on the Buslaev and Perelman methods \ci{BP1,BP2}: the symplectic
orthogonal decomposition of the dynamics near the solitary manifold,
the time decay for the linearized equation, etc.
Our problem differs from \cite{ikv} in the following aspects:

i) Speed of propagation for the   Schr\"odinger equation is
infinite, and the solitons
   exist only for the velocities $|v|<2m$.

ii) We consider nonspherically symmetric coupled function $\rho(x)$.
In this case we need additional arguments for
the absence of embedded
eigenvalues in the continuous spectrum.

iii) We also consider the coupling function $\rho(x)$ without compact
support. Respectively, for the proof of the time decay for the
linearized equation,  we use the Jensen-Kato results \ci{J, JK} and
the Agmon weighted norms \ci{Ag}.
\begin{remark}
The term $m^2$ in Schr\"odinger equation appears automatically in the
nonrelativistic limit of the Klein-Gordon equation, and traditionally is
removed by a gauge transformation. We keep the term to provide the existence
of the nonzero solitons.
\end{remark}


\setcounter{equation}{0}

\section{Main Results}


\subsection{Existence of Dynamics}

To formulate our results precisely, we need some definitions. We introduce
a suitable phase space for the Cauchy problem corresponding to (\ref{SS}) and (\ref{Hamil}).
Let $H^0=L^2$, and $H^1$ be the Sobolev space $H^1=\{\psi\in L^2:\,|\nabla\psi|\in L^2\}$
with the norm
$\Vert\psi\Vert_{H^1}=\Vert\nabla\psi\Vert_{L^2}+\Vert\psi\Vert_{L^2}$.
Let us introduce also the weighted Sobolev spaces $H^s_{\alpha}$, $s=0,1$,
$\alpha\in\R$ with the norms $\Vert\psi\Vert_{s,\alpha}:=\Vert(1+|x|)^{\alpha}\psi\Vert_{H^s}$.

\begin{definition}
  i) The phase space ${\cal E}$ is the real Hilbert space $H^1\oplus H^1\oplus {\R}^3\oplus {\R}^3$
  of states  $Y=(\psi_1 ,\psi_2 ,q,p)$ with the finite norm
  $$
  \Vert Y\Vert_{\cal E}=\Vert \psi_1 \Vert_{H^1} +
  \Vert\psi_2 \Vert_{H^1}+|q|+|p|.
  $$
  ii) ${\cal E}_{\alpha}$ is the space $H^1_{\alpha}\oplus H^1_{\alpha}\oplus {\R}^3\oplus {\R}^3$
  with the norm
  $$
  \Vert Y\Vert_{\alpha}=\Vert \,Y\Vert_{{\cal E}_{\alpha}}=
  \Vert \psi_1 \Vert_{1,\alpha} +\Vert\psi_2 \Vert_{1,\alpha}+|q|+|p|.
  $$
  iii) ${\cal E}^+$ is  space $H^2\oplus H^2\oplus {\R}^3\oplus {\R}^3$ with the  norm
  $$
  \Vert Y\Vert_{\cal E}^+=\Vert \psi_1 \Vert_{H^2} + \Vert\psi_2 \Vert_{H^2}+|q|+|p|.
  $$
\end{definition}
For $\psi_j\in L^2$ we have
\be\la{hbound}
  -\fr{1}{2m^2}\Vert\rho_j\Vert^2_{L^2}\le\frac{m^2}2\Vert\psi_j\Vert^2_{L^2}+
  \langle\psi_j,\rho_j(\cdot-q)\rangle\le\fr{m^2+1}{2}\Vert\psi_j\Vert^2_{L^2}+
  \fr{1}{2}\Vert\rho_j\Vert^2_{L^2}.
\ee
Therefore ${\cal E}$ is the space of finite energy states. The Hamilton functional
$\cal H$ is continuous on the space ${\cal E}$ and the lower bound in (\re{hbound})
implies that the energy (\re{Hamil}) is bounded from below.

The system (\re{SS}) reads as the Hamilton system
\be\la{ham}
  \dot Y=J{\cal D}{\cal H}(Y),\quad J:=\left(
  \begin{array}{cccc}
  0 & 1 & 0 & 0\\
  -1 & 0 & 0 & 0\\
  0 & 0 & 0 & 1\\
  0 & 0 & -1 & 0\\
  \end{array}
  \right),\quad Y=(\psi_1,\psi_2,q,p)\in{\cal E},
\ee
where ${\cal D}{\cal H}$ is the Fr\'echet derivative of the Hamilton functional (\re{Hamil}).
\begin{prop}\la{WPexistence}
  Let (\re{ro}) hold. Then\\
  (i) For every $Y_0\in {\cal E}$ the Cauchy problem (\re{WP2.1}) has a unique
  solution $Y(t)\in C(\R, {\cal E})$.\\
  (ii) For every $t\in\R$, the map $U(t): Y_0\mapsto Y(t)$ is continuous on ${\cal E}$.\\
  (iii) The  energy is conserved, i.e.
  \be\la{2.4}
     {\cal H}(Y(t))= {\cal H}(Y_0),\,\,\,\,\,t\in\R.
  \ee
\end{prop}
\begin{proof}
{\it Step i)}
Let us fix an arbitrary $b>0$ and prove (i)-(iii) for $Y_0\in{\cal E}$ such that
$\Vert Y_0\Vert_{\cal E}\le b$ and $|t|\le\ve=\ve(b)$ for some sufficiently small $\ve(b)>0$.
Let us rewrite the Cauchy problem \eqref{WP2.1} us
\be\la{WP2.2}
  \dot Y(t)=F_1(Y(t))+F_2(Y(t)),\quad t\in\R:\quad Y(0)=Y_0,
\ee
where $F_1: Y\mapsto ((-\Delta+m^2)\Psi_2,(\Delta-m^2)\Psi_1,0,0)$.
The Fourier transform provides the existence and uniqueness of solution
$Y_1(t)\in C(\R,{\cal E})$ to the linear problem \eqref{WP2.2} with $F_2=0$.
Let $U_1(t): Y_0\mapsto Y_1(t)$ be the corresponding strongly continuous
group of bounded linear operators on ${\cal E}$. Then \eqref{WP2.2} for
$Y(t)\in C(\R,{\cal E})$ is equivalent to
\be\la{Dugamel}
  Y(t)=U_1(t)Y_0+\int\limits_0^t~ds~U_1(t-s)F_2(Y(s)),
\ee
because $F_2(Y(\cdot))\in  C(\R,{\cal E})$ in this case.
The latter follows from a local Lipschitz continuity of the map $F_2$ in ${\cal E}$:
for each $b>0$ there exist a $\varkappa=\varkappa(b)>0$ such that for all $Y,Z\in{\cal E}$
with $\Vert Y\Vert_{\cal E},\;\Vert Z\Vert_{\cal E}\le b$,
$$
  \Vert F_2(Y)-F_2(Z)\Vert_{\cal E}\le\varkappa\Vert Y-Z\Vert_{\cal E}.
$$
Therefore, by the contraction mapping principle, equation \eqref{Dugamel}
has a unique local solution $Y(\cdot)\in C([-\ve,\ve],{\cal E})$ with
$\ve>0$ depending only on $b$.\\
{\it Step ii)}
We use now energy conservation to ensure the existence of a global solution and its continuity.
First consider $Y_0\in{\cal E}_c:= C_0^{\infty}\oplus  C_0^{\infty}\oplus {\R}^3\oplus {\R}^3$.
Then $Y(t)\in{\cal E^+}$ since $U_1(t)Y_0$, $F_2(Y(t))\in{\cal E^+}$ by \eqref{ro}.
The energy conservation law follows by \eqref{ham} and the the chain rule for the
Fr\'echet derivatives:
$$
  \ds\fr d{dt}{\cal H}(Y(t))=\langle D{\cal H}(Y(t)),\dot Y(t)\rangle=
  \langle D{\cal H}(Y(t)),J D{\cal H}(Y(t))\rangle=0,~~~~~~t\in\R
$$
since the operator $J$ is skew-symmetric by (\re{ham}), and
$D{\cal H}(Y(t))\in L^2\oplus L^2\oplus {\R}^3\oplus {\R}^3$ for  $Y(t)\in {\cal E^+}$.
The inequality \eqref{hbound} implies
$$
  {\cal H}\ge\fr 12\Vert\nabla\psi\Vert^2_{L^2}+\fr{m^2}4\Vert\psi\Vert^2_{L^2}+\frac 12|p|^2
  -\fr 1{m^2}\Vert\rho\Vert^2_{L^2}.
$$
Hence, by energy conservation, for $|t|\le\ve$
$$
  \fr 12\Vert\nabla\psi\Vert^2_{L^2}+\fr{m^2}4\Vert\psi\Vert^2_{L^2}+\frac 12|p|^2
  -\fr 1{m^2}\Vert\rho\Vert^2_{L^2}\le{\cal H}(Y(t))={\cal H}(Y_0).
$$
This implies  {\it a priori} estimate
\be\la{estimate}
  \Vert\psi\Vert_{H^1}+|p|\le B\quad {\rm for}\; |t|\le\ve,
\ee
with $B$ depending only on the norm $\Vert Y_0\Vert_{\cal E}$ of
the initial data and on $\Vert\rho\Vert_{L^2}$. An arbitrary
initial data $Y_0\in{\cal E}$ can be approximated by initial data
from ${\cal E}_c$. The corresponding solution exists due to
representation \eqref{Dugamel} by contraction
mapping principle, and then \eqref{estimate} follows by the limit transition.\\
{\it Step iii)}
Properties (i)-(iii) for arbitrary $t\in\R$ now follow from the same properties for small $|t|$
and from  {\it a priori} bound \eqref{estimate}.
\end{proof}

\subsection{Solitary Manifold and Main Result}

Let us compute the  solitons (\re{sosol}). The substitution to (\re{S}) gives the following
stationary equations
\be\la{stfch}
  \left.\begin{array}{l}
  -iv\cdot\na\psi_v(y)=(-\Delta+m^2)\psi_v(y)+\rho(y)
  \\
  p=v,\quad
  0=-\ds\int\big(\overline{\na\psi}_v(y)\rho(y)
  +\na\psi_v(y)\overline{\rho}(y)\big)\,dy
  \end{array}\right|
\ee
Then the first  equation implies
\be\la{Lpv}
  \Lambda\psi_v(y):=[-\Delta+m^2+iv\cdot\na]\psi_v(y)=-\rho(y),~~~~~~~~y\in\R^3.
\ee
For $|v|<2m$ the operator $\Lambda$ is an isomorphism  $H^4(\R^3)\to H^2(\R^3)$.
Hence  (\re{ro}) implies that
\be\la{Lpvy}
  \psi_v(y)=-\Lambda^{-1}\rho(y)\in H^4(\R^3).
\ee
If $v$ is given and $|v|<2m$, then $p_v$ can be found from the second equation of (\re{stfch}).

The function $\psi_v$ can be computed by the Fourier transform.
The soliton is given by the formula
\be\la{sol}
  \psi_v(x) = \ds -\frac 1{4\pi} \int\ds
  \frac {e^{-\sqrt{m^2-\fr{v^2}4}|x-y|}e^{i\fr v2(x-y)}\rho(y)d^3y}{|x-y|}.
\ee Further, in Appendix A, we prove that the last equation of of
(\re{stfch}) holds. Hence, the soliton solution (\re{sosol}) exists
and defined uniquely for any couple $(a,v)$ with $|v|<2m$. Let us
denote by  $V:=\{v\in\R^3:|v|<2m\}$, $\psi_{v1}=\Re\psi_v$, and
$\psi_{v2}=\Im\psi_v$.

\begin{definition}
  A soliton state is $S(\si):=(\psi_{v1}(x-b),\psi_{v2}(x-b),b,v)$, where
  $\si:=(b,v)$ with $b\in\R^3$ and $v\in V$.
\end{definition}
Obviously, the soliton solution admits the representation $S(\si(t))$, where
\be\la{sigma}
  \si(t)=(b(t),v(t))=(vt+a,v).
\ee
\begin{definition}
  A solitary manifold is the set ${\cal S}:=\{S(\si):\si\in\Sigma:=\R^3\times V\}$.
\end{definition}

The main result of our paper is the following theorem.
\begin{theorem}\la{main}
  Let \eqref{ro}, and the Wiener condition (\re{W}) hold. Let $\beta>3/2$ be the number from
  \eqref{ro}, and $Y(t)$ be the solution to the Cauchy problem  (\re{WP2.1}) with
  the initial state $Y_0$ which is sufficiently close to the solitary manifold:
  \be\la{close}
    p_0<2m,\quad d_0:={\rm dist}_{{\cal E}_\beta}(Y_0,{\cal S})\ll 1.
  \ee
  Then the asymptotics  hold for $t\to\pm\infty$,
  \be\la{nas}
    \dot q(t)=v_{\pm}+{\cal O}(|t|^{-2}),\quad q(t)=v_{\pm}t+a_{\pm}+{\cal O}(|t|^{-3/2});
  \ee
  \be\la{nasf}
    \psi(x,t)=\psi_{v\pm}(x-v_{\pm}t-a_{\pm})+W_0(t){\bp}_{\pm}+r_{\pm}(x,t)
  \ee
  with
  \be\la{rm}
     \Vert r_\pm(t)\Vert_{H^1}={\cal O}(|t|^{-1/2}).
  \ee
\end{theorem}
It suffices to  prove the asymptotics  (\re{nas}), (\re{nasf}) for $t\to+\infty$
since the system (\re{SS}) is time reversible.


\setcounter{equation}{0}

\section{
Symplectic Projection onto Solitary Manifold}

Let us identify the tangent space to ${\cal E}$, at every point,
with ${\cal E}$. Consider the symplectic form $\Om$ defined on
${\cal E}$ by $\Om=\ds\int d\psi_1(x)\wedge d\psi_2(x)\,dx+dq\wedge dp$, i.e.
\be\la{OmJ}
  \Om(Y_1,Y_2)=\langle Y_1,JY_2\rangle,\,\,\,Y_1,Y_2\in {\cal E},
\ee
where
$$
  \langle Y_1,Y_2\rangle:=\langle\psi_{11},\psi_{12}\rangle+
  \langle\psi_{21},\psi_{22}\rangle+q_1  q_2+p_1  p_2
$$
and $\langle\psi_{11},\psi_{12}\rangle=\ds\int\psi_{11}(x)\psi_{12}(x)dx$ etc.
It is clear that the form $\Om$ is non-degenerate, i.e.
$$
  \Om(Y_1,Y_2)=0\,\,~\mbox{\rm for every}~~ \,Y_2\in{\cal E} \,\,\Longrightarrow\,\, Y_1=0.
$$
\begin{definition}
  i) $Y_1\nmid Y_2$ means that $Y_1\in{\cal E}$, $Y_2\in{\cal E}$,
  and $Y_1$ is symplectic orthogonal to $Y_2$, i.e. $\Om(Y_1,Y_2)=0$.

  ii) A projection operator ${\bf P}:{\cal E}\to{\cal E}$ is called symplectic orthogonal if
  $Y_1\nmid Y_2$ for $Y_1\in\mbox{\rm Ker}\,{\bf P}$ and $Y_2\in \Im{\bf P}$.
\end{definition}


Let us consider the tangent space ${\cal T}_{S(\si)}{\cal S}$ to
the manifold ${\cal S}$ at a point $S(\si)$.
The vectors $\tau_j:=\pa_{\si_j}S(\si)$, where
$\pa_{\si_j}:=\pa_{b_j}$ and $\pa_{\si_{j+3}}:=\pa_{v_{j}}$ with $j=1,2,3$,
form a basis in ${\cal T}_{\si}{\cal S}$. In detail,
\be\la{inb}
\left.\begin{array}{rclrrrrcrcl}
  \tau_j=\tau_j(v)&:=&\pa_{b_j}S(\si)=
  (&\!\!\!\!-\pa_j\psi_{v1}(y)&\!\!\!\!,&\!\!\!\!-\pa_j\psi_{v2}(y)
   &\!\!\!\!,&\!\!e_j&\!\!\!\!,&\!\!0&\!\!\!\!)\\
   \tau_{j+3}=\tau_{j+3}(v)&:=&\pa_{v_j}S(\si)=
   (&\!\!\!\!\pa_{v_j}\psi_{v1}(y)&\!\!\!\!,&\!\!\!\!
   \pa_{v_j}\psi_{v2}(y)&\!\!\!\!,&\!\!0&\!\!\!\!,&\!\!
   e_j&\!\!\!\!)
\end{array}\right|\quad j=1,2,3,
\ee
where  $y:=x-b$ is the ``moving frame coordinate'', $e_1=(1,0,0)$ etc. Let us stress that
the functions $\tau_j$ will be considered always as the functions of $y$, not of $x$.

The formulas (\re{sol}) and the conditions (\re{ro}) imply that
\be\la{tana}
   \tau_j(v)\in{\cal E}_\al,\quad v\in V,\quad j=1,\dots,6,\quad \forall\al\le\beta.
\ee

\begin{lemma}\la{Ome}
  The matrix with the elements $\Om(\tau_l(v),\tau_j(v))$ is non-degenerate
  for any $v\in V$.
\end{lemma}
{\bf Proof } The elements are computed in Appendix B.
  As the result, the matrix $\Om(\tau_l,\tau_j)$ has the form
  \be\la{Omega}
     \Om(v):=(\Om(\tau_l,\tau_j))_{l,j=1,\dots,6}=\left(
     \begin{array}{ll}
         0 & \Om^+(v)\\
         -\Om^+(v) & 0
     \end{array}\right),
  \ee
  where the $3\times3$-matrix $\Om^+(v)$ equals
  \be\la{Wm}
     \Om^+(v)=K+E.
  \ee
  Here $K$ is a symmetric $3\times3$-matrix with the elements
  \be\la{alpha}
     K_{ij}=\int \frac{k_jk_l\Big((k^2+m^2)(|\hat\psi_{v1}|^2+|\hat\psi_{v2}|^2)
     +i(kv)(\hat\psi_{v1}\overline{\hat\psi}_{v2}-\hat\psi_{v2}\overline{\hat\psi}_{v1})
    \Big)~dk}{(k^2+m^2)^2-(kv)^2}
  \ee
  where the ``hat'' denotes the Fourier transform (cf. (\re{W})).
  The matrix $K$ is the integral of the symmetric nonnegative definite
  matrix $k\otimes k=(k_ik_j)$ with a nonnegative weight.
  (The last statement is true since $k^2+m^2>|(kv)|$ for $|v|<2m$,
  and $|\hat\psi_{v1}+i\hat\psi_{v2}|^2=|\hat\psi_{v1}|^2+|\hat\psi_{v2}|^2
  - i(\hat\psi_{v1}\overline{\hat\psi}_{v2}-\hat\psi_{v2}\overline{\hat\psi}_{v1})\ge 0.$)
  Hence, the matrix $K$ is also nonnegative definite.
  Since the unite matrix $E$ is positive definite, the matrix $\Om^+(v)$ is  symmetric and
  positive definite, hence non-degenerate. Then
  the matrix $\Om(\tau_l,\tau_j)$ also is non-degenerate.
\hfill $\bo$
Let us introduce the translations $T_a:(\psi_1(\cdot),\psi_2(\cdot),q,p)\mapsto
(\psi_1(\cdot-a),\psi_2(\cdot-a),q+a,p)$, $a\in\R^3$.
Note that the manifold ${\cal S}$ is invariant with respect to the translations.

\begin{definition}
   i) For any $\al\in\R$ and $\overline p<2m$ denote by
   ${\cal E}_\al(\overline p)=\{Y=(\psi_1,\psi_2,q,p)\in{\cal E}_\al:|p|
   \le\overline p\}$. We set ${\cal E}(\overline p):={\cal E}_0(\overline p)$.
\\
  ii) For any  $\overline v<2m$ denote by
  $\Si(\overline v)=\{\si=(b,v):b\in\R^3, |v|\le \overline v\}$.
\end{definition}
The next Lemma provide that in a small neighborhood of the soliton
manifold ${\cal S}$ a ``symplectic orthogonal projection''
onto ${\cal S}$ is well-defined. The  proof is similar to the proof of the Lemma 3.4 in \cite{ikv}.
\begin{lemma}\la{skewpro}
  Let (\re{ro}) hold, $\al\in\R$. Then
  \\
  i) there exists a neighborhood ${\cal O}_{\al}({\cal S})$ of ${\cal S}$ in ${\cal E}_\al$
  and a map ${\bf\Pi}:{\cal O}_{\al}({\cal S})\to{\cal S}$ such that ${\bf\Pi}$ is uniformly
  continuous in the metric of ${\cal E}_\al$
  on ${\cal O}_{\al}({\cal S})\cap {\cal E}_\al(\overline p)$ with $\overline p<2m$,
  \be\la{proj}
     {\bf\Pi} Y=Y~~\mbox{for}~~ Y\in{\cal S}, ~~~~~\mbox{and}~~~~~
     Y-S \nmid {\cal T}_S{\cal S},~~\mbox{where}~~S={\bf\Pi} Y.
  \ee
  ii) ${\cal O}_{\al}({\cal S})$ is invariant with respect to the translations $T_a$, and
  \[ {\bf\Pi} T_aY=T_a{\bf\Pi} Y,~~~~~\mbox{for}~~Y\in{\cal O}_{\al}({\cal S})
  ~~\mbox{and}~~a\in\R^3. \]
  iii) For any $\overline p<2m$ there exists a $\overline v<2m$ s.t.${\bf\Pi} Y=S(\si)$
  with  $\si\in \Si(\overline v)$ for
  $Y\in{\cal O}_{\al}({\cal S})\cap {\cal E}_\al(\overline p)$
  \\
  iv) For any $\overline v< 2m$ there exists an $r_\al(\overline v)>0$ s.t.
  $S(\si)+Z\in{\cal O}_{\al}({\cal S})$ if $\si\in\Si(\overline v)$
  and $\Vert Z\Vert_\al<r_\al(\overline v)$.
\end{lemma}

\medskip

\noindent
We will call ${\bf\Pi}$ a symplectic orthogonal projection onto ${\cal S}$.

\begin{cor}
  The condition (\re{close}) implies that $Y_0=S+Z_0$ where $S=S(\si_0)={\bf\Pi} Y_0$, and
  \be\la{closeZ}
     \Vert Z_0\Vert_\beta \ll 1.
  \ee
  \end{cor}
{\bf Proof}
   Lemma \re{skewpro} implies that ${\bf\Pi} Y_0=S$ is well defined for small $d_0>0$.
   Furthermore, the condition (\re{close}) means that there exists a point $S_1\in{\cal S}$
   such that $\Vert Y_0-S_1\Vert_\beta=d_0$. Hence,
   $Y_0,S_1\in {\cal O}_{\beta}({\cal S})\cap{\cal E}_\beta(\overline p)$
   with a $\overline p< 2m$ which does not depend on $d_0$ for sufficiently small $d_0$.
   On the other hand, ${\bf\Pi} S_1= S_1$, hence the uniform continuity of the map ${\bf\Pi}$
   implies that $\Vert S_1- S\Vert_\beta\to 0$ as $d_0\to 0$. Therefore, finally,
   $\Vert Z_0\Vert_\beta=\Vert Y_0- S \Vert_\beta \le \Vert Y_0- S_1 \Vert_\beta+
   \Vert S_1-S  \Vert_\beta\le d_0+o(1)\ll 1$ for small $d_0$.\bo


\setcounter{equation}{0}

\section{Linearization on the Solitary Manifold}

Let us consider a solution to the system (\re{SS}), and split it as the sum
\be\la{dec}
   Y(t)=S(\si(t))+Z(t),
\ee
where $\si(t)=(b(t),v(t))\in\Sigma$ is an arbitrary smooth function of $t\in\R$.
In detail, denote $Y=(\psi_1,\psi_2,q,p)$ and $Z=(\Psi_1,\Psi_2,Q,P)$.
Then (\re{dec}) means that
\be \la{add}\left.
  \begin{array}{rclrcl}
     \psi_1(x,t)&=&\psi_{v(t)1}(x-b(t))+\Psi_1(x-b(t),t),&q(t)&=&b(t)+Q(t)\\
     \psi_2(x,t)&=&\,\psi_{v(t)2}(x-b(t))+\Psi_2(x-b(t),t),&p(t)&=&v(t)+P(t)
  \end{array}\right|
\ee
Let us substitute (\re{add}) to (\re{SS}), and linearize the equations in $Z$.
Later we will choose $S(\si(t))={\bf\Pi} Y(t)$, i.e. $Z(t)$ is symplectic
orthogonal to ${\cal T}_{S(\si(t))}{\cal S}$.

Setting $y=x-b(t)$ which is the ``moving
frame coordinate'', we obtain from (\re{add}) and (\re{SS}) that
\be\la{addeq}\left.
  \begin{array}{rcl}
    \dot\psi_1&=&\dot v\cdot \na_v\psi_{v1}(y)-\dot b\cdot \na\psi_{v1}(y)+
    \dot\Psi_1(y,t)-\dot b\cdot \na\Psi_1(y,t)\\\\
    &=&-\Delta\psi_{v2}(y)+m^2\psi_{v2}(y)-\Delta\Psi_2(y,t)
    +m^2\Psi_2(y,t)+\rho_2(y-Q)\\\\
    \dot\psi_2&=&\dot v\cdot \na_v\psi_{v2}(y)-\dot b \cdot\na\psi_{v2}(y)+
    \dot\Psi_2(y,t)-\dot b\cdot \na\Psi_2(y,t)\\\\
    &=&\Delta\psi_{v1}(y)-m^2\psi_{v1}(y)+\Delta\Psi_1(y,t)
    -m^2\Psi_1(y,t)-\rho_1(y-Q)\\\\
    \dot q&=&\dot b+\dot Q=v+P\\\\
    \dot p&=&\dot v+\dot P=-\langle\na(\psi_{vj}(y)+
    \Psi_j(y,t)),\rho_j(y-Q)\rangle.
\end{array}\right|
\ee
Let us to extract linear terms in $Q$. First note that
$\rho_j(y-Q)=\rho_j(y)-Q \cdot\na\rho_j(y)+N_j(Q),\;j=1,2$, where
$N_j(Q)=\rho_j(y-Q)-\rho_j(y)+Q\cdot \na\rho_j(y)$. The condition \eqref{ro} implies
that for $N_j(Q)$ the bound holds,
\be\la{Nj}
  \Vert N_j(Q)\Vert_{0,\beta}\le C_\beta(\overline Q)Q^2,\;j=1,2,
\ee
uniformly in $|Q|\le\overline Q$ for any fixed $\overline Q$, where $\beta$ is the parameter
from Theorem \ref{main}. Using the equations (\re{stfch}), we obtain from (\re{addeq})
the following equations for the components of the vector $Z(t)$:
\be\la{Phi}
\left.\begin{array}{rcl}
\dot \Psi_1(y,t) &= &-\Delta\Psi_2(y,t)+m^2\Psi_2(y,t)+
\dot b\cdot \na\Psi_1(y,t)-Q \cdot\na\rho_2(y)\\\\
&+&(\dot b-v)\cdot \na\psi_{v1}(y)-\dot v\cdot \na_v\psi_{v1}(y)+N_2
\\\\
\dot \Psi_2(y,t) &=& \Delta\Psi_1(y,t)-m^2\Psi_1(y,t)+
\dot b\cdot \na\Psi_2(y,t)+Q \cdot\na\rho_1(y)\\\\
&+&(\dot b-v)\cdot \na\psi_{v2}(y)-\dot v\cdot\na_v\psi_{v2}(y)-N_1
\\\\
\dot Q(t) &=&P+(v-\dot b)
\\\\
\dot P(t) &=&\langle\Psi_j(y,t),\na\rho_j(y)\rangle+
\langle\na\psi_{vj}(y),Q \cdot\na\rho_j(y)\rangle-\dot v+N_4(v,Z)
\end{array}\right|
\ee
where $N_4(v,Z)=-\langle\na\psi_{vj},N_j(Q)\rangle+
\langle\na\Psi_j,Q\cdot \na\rho_j\rangle-\langle\na\Psi_j,N_j(Q)\rangle$.
Clearly, $N_4(v,Z)$ satisfies the following estimate
\be\la{N4}
  |N_4(v,Z)|\le C_\beta(\rho,\overline v,\overline Q)\Big[Q^2+
\Vert\Psi_1\Vert_{1,-\beta}|Q|+\Vert\Psi_2\Vert_{1,-\beta}|Q|
 \Big],
\ee
uniformly in $|v|\le \overline v$ and $|Q|\le \overline Q$
for any fixed  $\overline v<2m$.
We can write the equations (\re{Phi}) as
\be\la{lin}
\dot Z(t)=A(t)Z(t)+T(t)+N(t),\,\,\,t\in\R.
\ee
Here the operator $A(t)=A_{v,w}(t)$ depends on two parameters, $v=v(t)$, and
$w:=\dot b(t)$ and can be written in the form
\be\la{AA}
  A_{v,w}\left(
  \begin{array}{c}
     \Psi_1 \\ \Psi_2 \\ Q \\ P
\end{array}
\right)\\
=\left(
  \begin{array}{cccc}
     w \cdot\na & -(\Delta-m^2) & -\na\rho_2\cdot & 0 \\
     \Delta-m^2 & w\cdot \na & \na\rho_1\cdot & 0 \\
          0     &     0      &        0       & E \\
   \langle\cdot,\na\rho_1\rangle & \langle\cdot,\na\rho_2\rangle
   & \langle\na\psi_{vj},\cdot\na\rho_j\rangle & 0
  \end{array}
   \right)\left(
   \begin{array}{c}
      \Psi_1 \\ \Psi_2 \\ Q \\ P
  \end{array}
  \right).
\ee

Furthermore,   $T(t)=T_{v,w}(t)$ and $N(t)=N(t,\si,Z)$ in  (\re{lin}) stand for
\be\la{TN}
  T_{v,w}=\left(
  \begin{array}{c}
  (w-v)\cdot\na\psi_{v1}-\dot v\cdot\na_v\psi_{v1}\\
  (w-v)\cdot\na\psi_{v2}-\dot v\cdot\na_v\psi_{v2}\\
  v-w \\
  -\dot v
  \end{array}
  \right),\quad
  N(\si,Z)=\left(
  \begin{array}{c}
  N_2(Z) \\- N_1(Z) \\ 0 \\ N_4 (v,Z)
  \end{array}
  \right),
\ee
where $v=v(t)$, $w=w(t)$, $\si=\si(t)=(b(t),v(t))$, and $Z=Z(t)$.
The estimates (\re{Nj})  and (\re{N4}) imply that
\be\la{N14}
  \Vert N(\si,Z)\Vert_\beta\le C(\overline v, \overline Q)
  \Vert Z\Vert_{-\beta}^2,
\ee
uniformly in $\si\in\Si(\overline v)$ and $\Vert Z\Vert_{-\beta}\le r_{-\beta}(\overline v)$
for any fixed  $\overline v<2m$.
\begin{remark}\la{rT}
  {\rm
  i) The term $A(t)Z(t)$ in the right hand side of the equation  (\re{lin})
  is linear  in $Z(t)$, and $N(t)$ is a {\it high order term} in $Z(t)$.
  On the other hand, $T(t)$ is a zero order term which does not vanish at $Z(t)=0$
  since $S(\si(t))$ generally is not a soliton solution if (\re{sigma})
  does not hold (though $S(\si(t))$ belongs to the solitary manifold).
  \\
  ii) Formulas (\re{inb}) and (\re{TN}) imply:
  \be\la{Ttang}
     T(t)=-\sum\limits_{l=1}^3[(w-v)_l\tau_l+\dot v_l\tau_{l+3}]
  \ee
  and hence $T(t)\in {\cal T}_{S(\si(t))}{\cal S}$, $t\in\R$.
  This fact suggests an unstable character of the nonlinear dynamics
  {\it along the solitary manifold}.
  }
\end{remark}

 \setcounter{equation}{0}
\section{The Linearized Equation}

Here we collect some Hamiltonian and spectral properties of the generator (\re{AA}) of the
linearized equation. First, let us consider the linear equation
\be\la{line}
   \dot X(t)=A_{v,w}X(t),\quad t\in\R,\quad v\in V,\quad w\in \R^3.
\ee

\begin{lemma} \la{haml}(cf. \cite{ikv})
  i) For any $v\in V$ and $w\in \R^3$ the equation (\re{line}) can be written as the Hamilton system
  (cf. (\re{ham})),
  \be\la{lineh}
     \dot X(t)=
     JD{\cal H}_{v,w}(X(t)),~~~~~~~t\in\R,
  \ee
  where $D{\cal H}_{v,w}$ is the Fr\'echet derivative of the Hamilton functional
  \begin{multline}\la{H0}
    {\cal H}_{v,w}(X)=\fr12\int\Big[|\na\Psi_1|^2+m^2|\Psi_1|^2+
    |\na\Psi_2|^2+m^2|\Psi_2|^2\Big]dy+\int\Psi_2 w\cdot\na\Psi_1 dy\\
     +\int\rho_j(y)Q\cdot\na\Psi_j dy+
    \fr12P^2-\fr12\langle Q\cdot\na\psi_{vj}(y),Q\cdot\na\rho_j(y)\rangle,\quad
      X=(\Psi_1,\Psi_2,Q,P)\in {\cal E},
  \end{multline}
  \\
  ii) Energy conservation law holds for the solutions $X(t)\in C^1(\R,{\cal E}^+)$,
  \be\la{enec}
     {\cal H}_{v,w}(X(t))=\const,~~~~~t\in\R.
  \ee
  iii) The skew-symmetry relation holds,
  \be\la{com}
     \Omega(A_{v,w}X_1,X_2)=-\Omega(X_1,A_{v,w}X_2), ~~~~~~~~X_1,X_2\in {\cal E}.
  \ee
  i{\rm v})
  The operator $A_{v,w}$ acts on the tangent vectors $\tau_j(v)$ to the solitary
  manifold as follows,
  \begin{multline}\la{Atanform}
     A_{v,w}[\tau_j(v)]=(w-v)\cdot\na\tau_j(v),\,\,\,A_{v,w}[\tau_{j+3}(v)]=
     (w-v)\cdot\na\tau_{j+3}(v)+\tau_j(v),\;j=1,2,3.
  \end{multline}
\end{lemma}

We will apply Lemma \re{haml} mainly to the operator $A_{v,v}$ corresponding to $w=v$.
In that case the linearized equation has the following additional essential features.
\begin{lemma}\la{ceig}
  Let us assume that $w=v\in V$. Then\\
  i) The tangent vectors $\tau_j(v)$ with $j=1,2,3$ are eigenvectors,
  and $\tau_{j+3}(v)$ are root vectors of the
  operator $A_{v,v}$, corresponding to zero eigenvalue, i.e.
  \be\la{Atanformv}
    A_{v,v}[\tau_j(v)]=0,\,\,\,A_{v,v}[\tau_{j+3}(v)]=
    \tau_j(v),\,\,\,j=1,2,3.
  \ee
  ii) The Hamilton function (\re{H0}) is nonnegative definite since
  \be\la{H0vv}
     {\cal H}_{v,v}(X)=\ds\fr12\int|\Lambda^{1/2}(\Psi_1+i\Psi_2)
     -\Lambda^{-1/2}Q\cdot\na(\rho_1+i\rho_2)|^2dx+\fr12P^2\ge 0.
  \ee
  Here $\Lambda$ is the operator (\re{Lpv}) which is symmetric and nonnegative
  definite in $L^2(\R^3)$ for $|v|<2m$, and $\Lambda^{1/2}$ is the
  nonnegative definite square root defined in the Fourier representation.
\end{lemma}
\begin{proof}
  The first statement follows from (\re{Atanform}) with $w=v$.
  In order to prove ii) we rewrite the integral in (\re{H0vv})  as follows:
  \begin{multline}\la{H0vve}
    \fr 12\langle\Lambda^{1/2}(\Psi_1+i\Psi_2)-\Lambda^{-1/2}Q\cdot\na(\rho_1+i\rho_2),
    \Lambda^{1/2}(\Psi_1+i\Psi_2)-\Lambda^{-1/2}Q\cdot\na(\rho_1+i\rho_2)\rangle\\
    =\ds\fr 12\langle\Lambda(\Psi_1+i\Psi_2),\Psi_1+i\Psi_2\rangle-
    \langle\Psi_j,Q\cdot\na\rho_j\rangle
    +\ds\fr 12\langle\Lambda^{-1}Q\cdot\na(\rho_1+i\rho_2), Q\cdot\na(\rho_1+i\rho_2)\rangle
  \end{multline}
  since the operator $\Lambda^{1/2}$ is symmetric in $L^2(\R^3)$.
  Now all the terms of the expression (\re{H0vve}) can be identified with the
  corresponding terms in (\re{H0}) since
 \begin{multline*}
    \ds\fr 12\langle\Lambda(\Psi_1+i\Psi_2),\Psi_1+i\Psi_2\rangle=
    \ds\fr 12\langle[-\Delta+m^2+iv\cdot\na](\Psi_1+i\Psi_2),(\Psi_1+i\Psi_2)\rangle\\
    =\ds\fr 12\langle[-\Delta+m^2]\Psi_1,\Psi_1\rangle
    +\ds\fr 12\langle[-\Delta+m^2]\Psi_2,\Psi_2\rangle
    +\langle\Psi_2,v\cdot\nabla\Psi_1\rangle
 \end{multline*}
 and $\Lambda^{-1}(\rho_1+i\rho_2)=-(\psi_{v1}+i\psi_{v2})$ by (\re{Lpv}) and (\re{Lpvy}).
\end{proof}

\begin{remark}
{\rm
  For a soliton solution of the system(\re{SS}) we have $\dot b=v$, $\dot v=0$,
  and hence $T(t)\equiv 0$. Thus, the equation(\re{line}) is the linearization
  of the system (\re{SS}) on a soliton solution. In fact, we do not linearize
  (\re{SS}) on a soliton solution, but on a trajectory $S(\si(t))$ with $\si(t)$
  being nonlinear in $t$. We will show later that $T(t)$ is quadratic in $Z(t)$
  if we choose $S(\si(t))$ to be  the symplectic orthogonal projection of $Y(t)$.
  Then (\re{line})  is again the linearization of (\re{SS}).
}
\end{remark}


\setcounter{equation}{0}

\section{Symplectic Decomposition of the Dynamics}

Here we decompose the dynamics in two components: along the manifold ${\cal S}$ and in
transversal directions. The equation (\re{lin}) is obtained without any assumption on $\si(t)$
in (\re{dec}). We are going to choose $S(\si(t)):={\bf\Pi} Y(t)$, but then we need to know that
\be\la{YtO}
  Y(t)\in {\cal O}_{-\beta}({\cal S}),~~~~~t\in\R,
\ee
It is true for $t=0$ by our main assumption (\re{close}) with sufficiently small $d_0>0$.
Then  $S(\si(0))={\bf\Pi} Y(0)$ and $Z(0)=Y(0)-S(\si(0))$ are well defined.
We will prove below that (\re{YtO}) holds if $d_0$ is sufficiently small.
Let us choose an arbitrary $\overline v$ such that $|v(0)|<\overline v<2m$
and let $\delta=\overline v-|v(0)|$.
Denote by $r_{-\beta}(\overline v)$ the positive numbers from Lemma \re{skewpro} iv) which
corresponds to $\al=-\beta$. Then $S(\si)+Z\in {\cal O}_{-\beta}({\cal S})$ if $\si=(b,v)$
with $|v|<\overline v$ and $\Vert Z\Vert_{-\beta}<r_{-\beta}(\overline v)$. Note that
$\Vert Z(0)\Vert_{-\beta}<r_{-\beta}(\overline v)$ if $d_0$ is sufficiently small. Therefore,
$S(\si(t))={\bf\Pi}Y(t)$ and $Z(t)=Y(t)-S(\si(t))$ are well defined for $t\ge 0$ so
small that $|v|<\overline v$ and $\Vert Z(t)\Vert_{-\beta} < r_{-\beta}(\overline v)$. This
is formalized by the following standard definition.
\begin{definition}
  $t_*$ is the ``exit time'',
  \be\la{t*}
     t_*=\sup \{t>0: \Vert Z(s)\Vert_{-\beta} < r_{-\beta}(\overline v),~~
     |v(s)-v(0)|<\delta,~~0\le s\le t \}.
  \ee
\end{definition}

One of our main goals is to prove that $t_*=\infty$ if $d_0$ is sufficiently small.
This would follow if we show that
\be\la{Zt}
  \Vert Z(t)\Vert_{-\beta}<r_{-\beta}(\overline v)/2,~~|v(s)-v(0)|<\delta/2,~~~0\le t < t_*.
\ee
Note that
\be\la{Qind}
  |Q(t)|\le\overline Q:= r_{-\beta}(\overline v), ~~~~~0\le t< t_*.
\ee
Now $N(t)$ in (\re{lin}) satisfies, by (\re{N14}), the following estimate,
\be\la{Nest}
  \Vert N(t)\Vert_{\beta}\le C_\beta(\overline v)\Vert Z(t)\Vert^2_{-\beta},
  \,\,\,0\le t<t_*.
\ee


\subsection{Longitudinal Dynamics: Modulation Equations}

From now on we fix the decomposition $Y(t)=S(\si(t))+Z(t)$ for $0<t<t_*$
by setting $S(\si(t))={\bf\Pi} Y(t)$ which is equivalent to the
symplectic orthogonality condition of type (\re{proj}),
\be\la{ortZ}
  Z(t)\nmid{\cal T}_{S(\si(t))}{\cal S},\,\,\,0\le t<t_*.
\ee
This allows us to simplify drastically the asymptotic analysis of the dynamical equations
(\re{lin}) for the transversal component $Z(t)$. As the first step, we derive the longitudinal
dynamics, i.e. the ``modulation equations'' for the parameters $\si(t)$.
Let us derive a system of ordinary differential equations for the
vector $\si(t)$. For this purpose, let us write (\re{ortZ}) in the form
\be\la{orth}
  \Om(Z(t),\tau_j(t))=0,\,\,j=1,\dots,6, ~~~~~~~0\le t<t_*,
\ee
where the vectors $\tau_j(t)=\tau_j(\si(t))$ span the tangent space
${\cal T}_{S(\si(t))}{\cal S}$. Note that $\si(t)=(b(t),v(t))$, where
\be\la{sit}
  |v(t)|\le \overline v<2m,~~~~~~~~~0\le t<t_*,
\ee
by Lemma \re{skewpro} iii). It would be convenient for us to use some other parameters
$(c,v)$ instead of $\si=(b,v)$, where $c(t)= b(t)-\ds\int^t_0 v(\tau)d\tau$ and
\be\la{vw}
  \dot c(t)=\dot b(t)-v(t)=w(t)-v(t), ~~~~~~~~~0\le t<t_*.
\ee
We do not need an explicit form of the equations for $(c,v)$ but the
following statement, which can be proved similar to the Lemma 6.2 in \cite{ikv}.

\begin{lemma}\la{mod}
   Let $Y(t)$ be a solution to the Cauchy problem (\re{WP2.1}), and (\re{dec}),
   (\re{orth}) hold. Then $(c(t),v(t))$ satisfies the equation
   \be\la{parameq}
     \left(
     \begin{array}{l}
     \dot c(t) \\ \dot v(t)
     \end{array}
     \right)={\cal N}(\si(t),Z(t)), ~~~~~~~0\le t<t_*,
   \ee
  where
\be\la{NZ}
  {\cal N}(\si,Z)={\cal O}(\Vert Z\Vert^2_{-\beta})
\ee
uniformly in  $\si\in\Si(\overline v)$.
\end{lemma}


\subsection{Decay for the Transversal Dynamics}

In Section 11 we will show that our main Theorem \re{main} can be derived from
the following time decay of the transversal component $Z(t)$:
\begin{prop}\la{pdec}
   Let all conditions of Theorem \re{main} hold. Then $t_*=\infty$, and
   \be\la{Zdec}
     \Vert Z(t)\Vert_{-\beta}\le \ds\fr {C(\rho,\overline v,d_0)}{(1+|t|)^{3/2}},
     ~~~~~t\ge0.
   \ee
\end{prop}
We will derive (\re{Zdec}) in Sections 7-11 from our equation (\re{lin}) for the
transversal component $Z(t)$. This equation can be specified using Lemma \re{mod}.
Indeed, the lemma implies that
\be\la{Tta}
  \Vert T(t)\Vert_{\beta}\le C(\overline v)\Vert Z(t)\Vert^2_{-\beta},
  ~~~~~~~~~0\le t<t_*,
\ee
by (\re{TN})  since $w-v=\dot c$. Thus (\re{lin}) becomes the equation
\be\la{reduced}
  \dot Z(t)=A(t)Z(t)+\tilde N(t), ~~~~~~~~~0\le t<t_*,
\ee
where $A(t)=A_{v(t),w(t)}$, and $\tilde N(t):=T(t)+N(t)$ satisfies the estimate
\be\la{redN}
   \Vert\tilde  N(t)\Vert_{\beta}\le C\Vert Z(t)\Vert^2_{-\beta},~~~~
   ~~~~~~~~~0\le t<t_*.
\ee
In all remaining part of our paper we will analyze mainly the {\bf basic equation}
(\re{reduced}) to establish the decay (\re{Zdec}). We are going to derive the decay using
the bound (\re{redN}) and the orthogonality condition  (\re{ortZ}).

First, we reduce the problem to the analysis of the {\it frozen} linear equation,
\be\la{Avv}
  \dot X(t)=A_1X(t), ~~t\in\R,
\ee
where $A_1$ is the operator $A_{v_1,v_1}$
defined by (\re{AA})
with $v_1=v(t_1)$ and a fixed $t_1\in[0,t_*)$. Then  we can apply wellknown
methods of scattering theory and then estimate the error by the method of majorants.

Note, that even for the frozen equation (\re{Avv}), the decay of type  (\re{Zdec})
for all solutions does not hold without  the orthogonality condition of type (\re{ortZ}).
Namely, by  (\re{Atanformv}) the equation (\re{Avv}) admits the {\it secular solutions}
\be\la{secs}
   X(t)=\sum_1^3 C_{j}\tau_j(v)+\sum_1^3 D_j[\tau_j(v)t+\tau_{j+3}(v)]
\ee
which arise also by differentiation of the soliton (\re{sosol}) in the parameters
$a$ and $v$ in the moving coordinate $y=x-v_1t$. Hence, we have to take into account
the orthogonality condition  (\re{ortZ}) in order to avoid the secular solutions.
For this purpose we will apply the corresponding symplectic orthogonal projection
which kills the ``runaway solutions''  (\re{secs}).
\begin{remark}
{\rm
  The solution (\re{secs}) lies in the tangent space  ${\cal T}_{S(\si_1)}{\cal S}$
  with $\si_1=(b_1,v_1)$ (for an arbitrary $b_1\in\R$) that suggests an unstable
  character of the nonlinear dynamics {\it along the solitary manifold}
  (cf. Remark \re{rT} iii)).
}
\end{remark}

\begin{definition}
  i) For $v\in V$, denote by ${\bf\Pi}_v$ the symplectic orthogonal projection
  of ${\cal E}$ onto the tangent space ${\cal T}_{S(\si)}{\cal S}$, and
  ${\bf P}_v={\bf I}-{\bf\Pi}_v$.\\
  ii) Denote by ${\cal Z}_v={\bf P}_v{\cal E}$ the space symplectic orthogonal to
  ${\cal T}_{S(\si)}{\cal S}$ with $\si=(b,v)$ (for an arbitrary $b\in\R$).
\end{definition}
Note that by the linearity,
\be\la{Piv}
  {\bf\Pi}_vZ=\sum{\bf\Pi}_{jl}(v)
  \tau_j(v)\Om(\tau_l(v),Z),~~~~~~~~~~Z\in{\cal E},
\ee
with some smooth coefficients ${\bf\Pi}_{jl}(v)$. Hence, the projector ${\bf\Pi}_v$,
in the variable $y=x-b$, does not depend on $b$, and this explains the choice
of the subindex in ${\bf\Pi}_v$ and ${\bf P}_v$.

Now we have the symplectic orthogonal decomposition
\be\la{sod}
  {\cal E}={\cal T}_{S(\si)}{\cal S}+{\cal Z}_v,~~~~~~~\si=(b,v),
\ee
and the symplectic orthogonality  (\re{ortZ}) can be written in the following
equivalent forms,
\be\la{PZ}
  {\bf\Pi}_{v(t)} Z(t)=0,~~~~{\bf P}_{v(t)}Z(t)= Z(t),~~~~~~~~~0\le t<t_*.
\ee

\begin{remark}\la{rZ}
{\rm
  The tangent space ${\cal T}_{S(\si)}{\cal S}$ is invariant under the operator
  $A_{v,v}$ by Lemma \re{ceig} i), hence the space  ${\cal Z}_v$ is also invariant
  by (\re{com}): $A_{v,v}Z\in {\cal Z}_v$ for {\it sufficiently smooth}  $Z\in{\cal Z}_v$.
}
\end{remark}

In Sections 12-18 we will prove the following proposition which will be one of the main
ingredients for proving \eqref{Zdec}. Let us consider the Cauchy problem for the equation \eqref{Avv}
with $A=A_{v,v}$ for a fixed $v\in V$. Recall that the parameter $\beta>3/2$ is also fixed.
\begin{prop}\la{lindecay}
  Let the conditions \eqref{ro} and \eqref {W} hold, $|v|\le\overline v<2m$, and $X_0\in{\cal E}$.
  Then\\
  i) The equation (\re{Avv}), with $A=A_{v,v}$, admits the unique solution
  $e^{At}X_0:=X(t)\in C(\R, {\cal E})$ with the initial condition $X(0)=X_0$.\\
  ii) For $X_0\in{\cal Z}_v\cap{\cal E_{\beta}}$, the solution
  $X(t)$ has the following decay,
  \be\la{frozenest}
     \Vert e^{At}X_0\Vert_{-\beta}\le
     \fr{C_{\beta}(\overline v)}{(1+|t|)^{3/2}}\Vert X_0\Vert_{\beta},\quad t\in\R.
  \ee
\end{prop}

\setcounter{equation}{0}

\section{Frozen Transversal Dynamics}

Now  let us fix an arbitrary $t_1\in [0,t_*)$, and
rewrite the equation (\re{reduced}) in a ``frozen form''
\be\la{froz}
   \dot Z(t)=A_1Z(t)+(A(t)-A_1)Z(t)+\tilde N(t),\,\,\,~~~~0\le t<t_*,
\ee
where $A_1=A_{v(t_1),v(t_1)}$ and
\begin{equation*}
   A(t)-A_1 =\left(
  \begin{array}{cccc}
   [w(t)-v(t_1)]\cdot \na & 0 & 0 & 0 \\
   0 & [w(t)-v(t_1)]\cdot\na & 0 & 0 \\
   0 & 0 & 0 & 0 \\
   0 & 0 & \langle\na(\psi_{v(t)j}-\psi_{v(t_1)j}),\na\rho_j\rangle & 0
  \end{array}
  \right).
\end{equation*}
The next trick is important since it allows us to kill the ``bad terms''
 $[w(t)\!-\!v(t_1)]\cdot \na$ in the operator $A(t)-A_1$.
\begin{definition}\la{d71}
  Let us change the  variables $(y,t)\mapsto (y_1,t)=(y+d_1(t),t)$, where
  \be\la{dd1}
     d_1(t):=\int_{t_1}^t(w(s)-v(t_1))ds, ~~~~0\le t\le t_1.
  \ee
\end{definition}
Next define
\be\la{Z1}
   Z_1(t):=(\Psi_1(y_1-d_1(t),t),\Psi_2(y_1-d_1(t),t),Q(t),P(t)).
\ee
Then we obtain the final form of the ``frozen equation'' for the transversal dynamics
\be\la{redy1}
  \dot Z_1(t)=A_1Z_1(t)+B_1(t)Z_1(t)+\tilde N_1(t),\,\,\,0\le t\le t_1,
\ee
where $\tilde N_1(t)=\tilde N(t)$ expressed in terms of $y=y_1-d_1(t)$,  and
$$
  B_1(t)=\left(
  \begin{array}{cccc}
  0 & 0 & 0 & 0 \\
  0 & 0 & 0 & 0 \\
  0 & 0 & 0 & 0 \\
  0 & 0 & \langle\na(\psi_{v(t)j}\!-\!\psi_{v(t_1)j}),\na\rho_j\rangle & 0
  \end{array}
  \right).
$$
\begin{lemma}\la{dest} (see \cite{ikv})
For $(\Psi_1,\Psi_2,Q,P)\in{\cal E}_{\alpha}$ with any $\alpha\le\beta$
the following estimate holds:
\be\la{shiftest}
\Vert(\Psi_1(y_1-d_1),\Psi_2(y_1-d_1),Q,P)\Vert_{\alpha}\le
\Vert(\Psi_1,\Psi_2,Q,P)\Vert_{\alpha}(1+|d_1|)^{|\alpha|}~,\,\,\,~~~~~~d_1\in\R^3.
\ee
\end{lemma}
\begin{cor}\la{destc}
 The following  bounds hold  for  $0\le t\le t_1$
   \be\la{dest1}
      \Vert\tilde N_1(t)\Vert_{\beta}\le\Vert Z(t)\Vert^2_{-\beta}(1+|d_1(t)|)^{\beta}~,\quad
      \Vert B_1(t)Z_1(t)\Vert_{\beta}\le C\Vert Z(t)\Vert_{-\beta}\int_t^{t_1}
      \Vert Z(\tau)\Vert^2_{-\beta}d\tau~.
   \ee
\end{cor}


\setcounter{equation}{0}

\section{Integral Inequality}

The equation (\re{redy1}) can be written in the integral form:
\be\la{Z1duh}
  Z_1(t)=e^{A_1t}Z_1(0)+\int_0^te^{A_1(t-s)}[B_1Z_1(s)+\tilde N_1(s)]ds,\quad 0\le t\le t_1.
\ee
Now we apply the symplectic orthogonal projection ${\bf P}_1:={\bf P}_{v(t_1)}$ to both sides
of \eqref{Z1duh}. The space ${\cal Z}_1:={\bf P}_1{\cal E}$ is invariant with respect to
$e^{A_1t}$ by Proposition \re{lindecay} ii) (cf. also Remark \re{rZ}). Therefore ${\bf P}_1$
commutes with the group $e^{A_1t}$ and applying (\re{frozenest}) we obtain that
$$
  \Vert {\bf P}_1Z_1(t)\Vert_{-\beta}
  \le C\fr{\Vert {\bf P}_1Z_1(0)\Vert_{\beta}}{(1+t)^{3/2}}
  +C\int_0^t\fr{\Vert{\bf P}_1[B_1Z_1(s)+\tilde N_1(s)]\Vert_{\beta}~ds}{(1+|t-s|)^{3/2}}.
$$
The operator ${\bf P}_1={\bf I}-{\bf\Pi}_1$ is continuous in ${\cal E}_\beta$ by
(\re{Piv}). Hence, using \re{dest1} we obtain that
\begin{multline}\la{duhest}
   \Vert {\bf P}_1Z_1(t)\Vert_{-\beta}\le\fr{C(\overline d_1(0))}{(1+t)^{3/2}}\Vert Z(0)\Vert_{\beta}\\
   +C(\overline d_1(t))\int_0^t\!\!\fr1{(1+|t-s|)^{3/2}}\left[\Vert Z(s)\Vert_{-\beta}
   \int_s^{t_1}\Vert Z(\tau)\Vert^2_{-\beta}d\tau+\Vert Z(s)\Vert^2_{-\beta}\right]ds,
   \quad 0\le t\le t_1.
 \end{multline}
where $\overline d_1(t):=\sup_{0\le s\le t} |d_1(s)| $.
\begin{definition} $t_{*}'$ is the exit time
  \be\la{t*'}
    t_*'=\sup \{t\in[0,t_*): \overline d_1(s)\le 1,~~0\le s\le t\}.
  \ee
\end{definition}
Now (\re{duhest}) implies that for $t_1<t_*'$
\begin{multline}\la{duhestri}
   \Vert {\bf P}_1Z_1(t)\Vert_{-\beta}\le\fr{C}{(1+t)^{3/2}}\Vert Z(0)\Vert_{\beta}\\
   +C_1\int_0^t\fr1{(1+|t-s|)^{3/2}}\left[\Vert Z(s)\Vert_{-\beta}
   \int_s^{t_1}\Vert Z(\tau)\Vert^2_{-\beta}d\tau+
   \Vert Z(s)\Vert^2_{-\beta}\right]ds,\quad 0\le t\le t_1.
\end{multline}

\setcounter{equation}{0}

\section{Symplectic Orthogonality}
Finally, we are going to change ${\bf P}_1Z_1(t)$ by $Z(t)$ in the left hand side of
(\re{duhestri}). We will prove that it is possible using again that $d_0\ll 1$ in (\re{close}).
For the justification we reduce further the exit time. First, we introduce the ``majorant''
\be\la{maj}
  m(t):=
  \sup_{s\in[0,t]}(1+s)^{3/2}\Vert Z(s)\Vert_{-\beta}~,~~~~~~~~~t\in [0,t_*).
\ee
Let us denote by $\ve$ a fixed positive number which we
will specify below.
\begin{definition} $t_{*}''$ is the exit time
  \be\la{t*''}
  t_*''=\sup \{t\in[0,t_*'):
  m(s)\le \ve,~~0\le s\le t\}.
\ee
\end{definition}

The following important bound  (\re{Z1P1est}) allows us to change the norm of
${\bf P}_1Z_1(t)$ in the left hand side of (\re{duhestri}) by the norm of $Z(t)$.

\begin{lemma}\la{Z1P1Z1}(cf.\cite{ikv})
  For sufficiently small $\ve>0$, we have for  $t_1<t_*''$
  \be\la{Z1P1est}
     \Vert Z(t)\Vert_{-\beta}\le C\Vert {\bf P}_1Z_1(t)\Vert_{-\beta},
     ~~~~~~~~0\le t \le t_1,
  \ee
  where $C$ depends only on $\rho$ and $\overline v$.
\end{lemma}
\begin{proof}
 Since $|d_1(t)|\le 1$ for $t\le t_1<t_*''<t_*'$ then by Lemma \re{dest} it suffices
 to prove that
 \be\la{Z1P1ests}
   \Vert Z_1(t)\Vert_{-\beta}\le 2\Vert{\bf P}_1Z_1(t)\Vert_{-\beta},\quad 0\le t\le t_1.
 \ee
 Recall that ${\bf P}_1Z_1(t)=Z_1(t)-{\bf\Pi}_{v(t_1)}Z_1(t)$.
 Then estimate \eqref{Z1P1ests} will follow from
 \be\la{Z1P1estf}
 \Vert{\bf\Pi}_{v(t_1)}Z_1(t)\Vert_{-\beta}
 \le\fr 12\Vert Z_1(t)\Vert_{-\beta},\quad 0\le t\le t_1.
 \ee
 Symplectic orthogonality \eqref{PZ} implies
 \be\la{PZ0r}
   {\bf\Pi}_{v(t),1}Z_1(t)=0,~~~~t\in[0,t_1],
 \ee
 where ${\bf\Pi}_{v(t),1}Z_1(t)$ is ${\bf\Pi}_{v(t)}Z(t)$ expressed in terms of
 the variable $y_1=y+d_1(t)$. Hence, (\re{Z1P1estf}) follows from (\re{PZ0r}) if the difference
 ${\bf\Pi}_{v(t_1)}-{\bf\Pi}_{v(t),1}$ is small uniformly in $t$, i.e.
 \be\la{difs}
   \Vert{\bf\Pi}_{v(t_1)}-{\bf\Pi}_{v(t),1}\Vert<1/2,~~~~~~~0\le t\le t_1.
 \ee
 It remains to justify (\ref{difs}) for small enough $\varepsilon>0.$
 Formula (\ref{Piv}) implies the following relation
 \be\la{Piv1}
  {\bf\Pi}_{v(t),1}Z_{1}(t)=\sum{\bf\Pi}_{jl}(v(t))\tau_{j,1}(v(t))
  \Omega (\tau_{l,1}(v(t)),Z_1(t)),
 \ee
 where $\tau_{j,1}(v(t))$ are the vectors $\tau_j(v(t))$ expressed in the variables $y_{1}$.
 Since $|d_1(t)|\le 1$ and $\nabla\tau_j$
 are smooth and fast decaying at infinity functions, Lemma \re{dest} implies that
 \be\la{011}
   \Vert \tau_{j,1}(v(t))-\tau_j(v(t))\Vert_\beta\le C|d_1(t)|^\beta,\quad 0\le t\le t_1
 \ee
 for all $j=1,2,\dots,6$. Furthermore,
 \[   \tau_j(v(t))-\tau_j(v(t_1)) =\int_t^{t_1}\dot v(s)\cdot\nabla_v\tau_j(v(s))ds,    \]
 and therefore
 \be\la{012}
   \Vert \tau_j(v(t))-\tau_j(v(t_1))\Vert_\beta
   \le C\int_t^{t_1}|\dot v(s)|ds,\quad 0\le t\le t_1.
 \ee
 Similarly,
 \be\la{013}
   {\bf\Pi}_{jl}(v(t))-{\bf\Pi}_{jl}(v(t_1))|= |\int_t^{t_1}\dot v(s)\cdot\nabla_v
   {\bf\Pi}_{jl}(v(s))ds|\le C\int_t^{t_1}|\dot v(s)|ds,~~~~0\leq t\le t_1,
 \ee
 since $|\nabla_v{\bf\Pi}_{jl}(v(s))|$ is uniformly bounded by (\re{sit}). Hence, the bounds
 (\ref{difs}) will follow from (\ref{Piv}), (\ref{Piv1}) and (\ref{011})-(\ref{013}) if we
 establish that $|d_1(t)|$ and the integral in the right hand side of (\ref{012}) can be made
 as small as we please by choosing $\ve >0$ small enough.

 To estimate $d_1(t)$, we note that
 \be\la{wen}
   w(s)-v(t_1)=w(s)-v(s)+v(s)-v(t_1)= \dot c(s)+\int_s^{t_1}\dot v(\tau)d\tau
 \ee
 by (\re{vw}). Hence, (\ref {dd1}), Lemma \ref{mod} and the definition (\re{maj}) imply that
 \begin{multline}\label{d1est}
   |d_1(t)| =|\int_{t_{1}}^t(w(s)-v(t_1))ds|\le
   \int_t^{t_{1}}\left( |\dot c(s)|+\int_s^{t_1}|\dot v(\tau)|d\tau\right)ds\\\\
   \le  Cm^2(t_1)\int_t^{t_1}\left(\frac 1{(1+s)^3}
   +\int_s^{t_1}\frac{d\tau}{(1+\tau)^3}\right) ds\le Cm^2(t_1)\le C\ve^2,\quad 0\le t\le t_1
 \end{multline}
 since $t_1<t_{\ast}^{\prime\prime}$. Similarly,
 \be\la{tvjest}
   \int_t^{t_1}|\dot v(s)|ds\le Cm^2(t_1)\int_t^{t_{1}}\frac{ds}{(1+s)^3}\le C\ve^2,\quad
   0\le t\le t_1.
\ee
\end{proof}

\setcounter{equation}{0}

\section{Decay of Transversal Component}
Here we prove Proposition \re{pdec}.
\\
{\it Step i)} We fix $\ve>0$ and $t''_*=t''_*(\ve)$ for which Lemma \re{Z1P1Z1} holds.
Then the bound of type (\re{duhestri}) holds with
$\Vert {\bf P}_1Z_1(t)\Vert_{-\beta}$ in the left hand side replaced by
$\Vert Z(t)\Vert_{-\beta}$~:
\begin{multline}\la{duhestrih}
   \Vert Z(t)\Vert_{-\beta}\le\fr{C}
   {(1+t)^{3/2}}\Vert Z(0)\Vert_{\beta}\\
   +C\int_0^t\fr1{(1+|t-s|)^{3/2}}\left[\Vert Z(s)\Vert_{-\beta}
   \int_s^{t_1}\Vert Z(\tau)\Vert^2_{-\beta}d\tau+
   \Vert Z(s)\Vert^2_{-\beta}\right]ds,\quad
   0\le t\le t_1
\end{multline}
for $t_1<t_*'$. This implies an integral inequality for the majorant $m(t)$
defined in \eqref{maj}.
Namely, multiplying both sides of (\re{duhestrih}) by $(1+t)^{3/2}$,
and taking the supremum in $t\in[0,t_1]$, we get
\begin{multline*}
\!\!\!\!  m(t_1)\le C\Vert Z(0)\Vert_{\beta}+C\sup_{t\in[0,t_1]}\ds
  \int_0^t\fr{(1+t)^{3/2}}{(1+|t-s|)^{3/2}}\left[\fr{m(s)}{(1+s)^{3/2}}
  \int_s^{t_1}\fr{m^2(\tau)d\tau}{(1+\tau)^{3}}+\fr{m^2(s)}{(1+s)^{3}}\right]ds
\end{multline*}
for $t_1\le t_*''$. Taking into account that $m(t)$ is a monotone
increasing function, we get
\be\la{mest}
   m(t_1)\le C\Vert Z(0)\Vert_{\beta}+C[m^3(t_1)+m^2(t_1)]I(t_1),
   ~~~~~~~~~~~~~t_1\le t_*''.
   \ee
where
\begin{multline*}
  I(t_1)=\sup_{t\in[0,t_1]}
  \int_0^{t}\fr{(1+t)^{3/2}}{(1+|t-s|)^{3/2}}\left[\fr1{(1+s)^{3/2}}
  \int_s^{t_1}\fr{d\tau}{(1+\tau)^{3}}+\fr1{(1+s)^3}\right]ds\le \overline I<\infty.
\end{multline*}
Therefore, (\re{mest}) becomes
\be\la{m1est}
  m(t_1)\le C\Vert Z(0)\Vert_{\beta}+C\overline I[m^3(t_1)+m^2(t_1)],~~~~ t_1<t_*''.
\ee
This inequality implies that $m(t_1)$ is bounded for $t_1<t_*''$, and moreover,
\be\la{m2est}
  m(t_1)\le C_1\Vert Z(0)\Vert_{\beta},~~~~~~~~~t_1<t_*''\,,
\ee
since $m(0)=\Vert Z(0)\Vert_{\beta}$ is sufficiently small by (\re{closeZ}).
\\
{\it Step ii)} The constant $C_1$ in the estimate (\re{m2est}) does not depend on
$t_*$, $t_*'$ and $t_*''$ by Lemma \re{Z1P1Z1}. We choose $d_0$ in (\re{close})
so small that $\Vert Z(0)\Vert_{\beta}<\ve/(2C_1)$. It is possible due to (\re{closeZ}).
Then the estimate (\re{m2est}) implies that $t''_*=t'_*$ and therefore
(\re{m2est}) holds for all $t_1<t'_*$. Then the bound (\re{d1est}) holds
for all $t<t_*'$. We choose $\ve$ so small that the right hand side in (\re{d1est})
does not exceed one. Then $t'_*=t_*$. Therefore, (\re{m2est}) holds for all $t_1<t_*$,
hence the first inequality in (\re{Zt}) also holds if $\Vert Z(0)\Vert_{\beta}$ is
sufficiently small by \eqref{maj} and \eqref{tvjest}.
Finally, this implies that $t_*=\infty$, hence also $t''_*=t'_*=\infty$ and
(\re{m2est}) holds for all $t_1>0$ if $d_0$ is small enough.
 \bo

\section{Soliton Asymptotics}
\setcounter{equation}{0}

Here we prove our main Theorem \re{main} under the assumption that the decay (\re{Zdec}) holds.
First we will prove the asymptotics (\re{qq}) for the vector components, and afterwards
the asymptotics (\re{S}) for the fields.
\\
{\bf Asymptotics for the vector components}
From (\re{addeq}) we have $\dot q=\dot b+\dot Q$, and from (\re{reduced}), (\re{redN}), (\re{AA})
it follows that $\dot Q=P+{\cal O}(\Vert Z\Vert^2_{-\beta})$. Thus,
\be\la{dq}
   \dot q=\dot b+\dot Q=v(t)+\dot c(t)+P(t)+{\cal O}(\Vert Z\Vert^2_{-\beta}).
\ee
The equation (\re{parameq}) and the estimates (\re{NZ}), (\re{Zdec}) imply that
\be\la{bv}
  |\dot c(t)|+|\dot v(t)|\le \ds\fr {C_1(\rho,\overline v,d_0)}{(1+t)^{3}},
  ~~~~~~t\ge0.
\ee
Therefore, $c(t)=c_+ +{\cal O}(t^{-2})$ and $v(t)=v_+ +{\cal O}(t^{-2})$, $t\to\infty$. Since
$|P|\le\Vert Z\Vert_{-\beta}$, the estimate (\re{Zdec}), and (\re{bv}), (\re{dq}) imply that
\be\la{qbQ}
  \dot q(t)=v_++{\cal O}(t^{-3/2}).
\ee
Similarly,
\be\la{bt}
  b(t)=c(t)+\ds\int_0^tv(s)ds=v_+t+a_++{\cal O}(t^{-1}),
\ee
hence the second part of (\re{qq}) follows:
\be\la{qbQ2}
  q(t)=b(t)+Q(t)=v_+t+a_++{\cal O}(t^{-1}),
\ee
since $Q(t)={\cal O}(t^{-3/2})$ by  (\re{Zdec}).
\\
{\bf Asymptotics for the fields}
We apply the approach developed in \ci{IKSs}, see also \ci{ikv}. For the field part of
the solution, $\psi(x,t)=\psi_1(x,t)+i\psi_2(x,t)$ let us define the accompanying soliton
field as $\psi_{\rm v(t)}(x-q(t))$, where we define now ${\rm v}(t)=\dot q(t)$, cf. (\re{dq}).
Then for the difference $z(x,t)=\psi(x,t)-\psi_{\rm v(t)}(x-q(t))$ we obtain the equation
$$
  i\dot z(x,t)=(-\Delta+m^2)z(x,t)-i\dot{\rm v}\cdot\na_{\rm v}\psi_{{\rm v}(t)}(x-q(t)).
$$
Then
\be\la{eqacc}
z(t)=W_0(t)z(0)-\int_0^tW_0(t-s)[i\dot{\rm v}(s)\cdot\na_{\rm v}\psi_{{\rm v}(s)}(\cdot-q(s))]ds.
\ee
To obtain the asymptotics (\re{nasf}) it suffices to prove that $z(t)=W_0(t){\bp}_++r_+(t)$ with
some ${\bp}_+\in H^1$ and $\Vert r_+(t)\Vert_{H^1}={\cal O}(t^{-1/2})$.
This is equivalent to
\be\la{Sme}
  W_0(-t)z(t)={\bp}_++r_+'(t),
\ee
where $\Vert r_+'(t)\Vert_{H^1}={\cal O}(t^{-1/2})$ since $W^0(t)$ is a unitary group in the
Sobolev space ${\cal F}$ by the energy conservation for the free Schr\"odinger equation.
Finally, (\re{Sme}) holds since (\re{eqacc}) implies that
$$
  W_0(-t)z(t)=z(0)-\int_0^t W_0(-s)f(s)ds,\quad
  f(s)=i\dot{\rm v}(s)\cdot\na_{\rm v}\psi_{{\rm v}(s)}(\cdot-q(s)),
$$
where the integral  in the right hand side of  converges in the Hilbert space ${\cal F}$
with the rate ${\cal O}(t^{-1/2})$. The latter holds since
$\Vert  W_0(-s)f(s)\Vert_{H^1} ={\cal O}(s^{-3/2})$ by the unitarity of $W_0(-s)$ and
the decay rate $\Vert f(s)\Vert_{H^1} ={\cal O}(s^{-3/2})$. Let us prove this rate of decay.
It suffices to prove that $|\dot {\rm v}(s)|={\cal O}(s^{-3/2})$,
or equivalently $|\dot p(s)|={\cal O}(s^{-3/2})$. Substitute (\re{add}) to the last
equation of (\re{SS}) and obtain
\begin{multline*}
\!\!\!\dot p(t)=\int\left[\psi_{v(t)j}(x-b(t))+\Psi_j(x-b(t),t)\right]\na\rho_j(x-b(t)-Q(t))dx
=\int\psi_{v(t)j}(y)\na\rho_j(y)dy\\
+\int\psi_{v(t)j}(y)\left[\na\rho_j(y-Q(t))-\na\rho_j(y)\right]dy
+\int\Psi_j(y,t)\na\rho_j(y-Q(t))dy.
\end{multline*}
The first integral in the right hand side is zero by the stationary equations (\re{stfch}).
The second integral is ${\cal O}(t^{-3/2})$, since $Q(t)={\cal O}(t^{-3/2})$, and by the conditions
(\re{ro}) on $\rho$. Finally, the third integral is ${\cal O}(t^{-3/2})$
by the estimate (\re{Zdec}). The proof is complete.\bo

 \setcounter{equation}{0}

\section{Decay for the Linearized Dynamics}
In remaining sections we prove Proposition \re{lindecay} in order to complete the proof of the
main result (Theorem \re{main}). Here we discuss our general strategy of the proof of the
Proposition. We apply the Fourier-Laplace transform
\be\la{FL}
  \tilde X(\lambda)=\int_0^\infty e^{-\lambda t}X(t)dt,~~~~~~~\Re\lambda>0
\ee
to (\re{Avv}).  According to Proposition \re{lindecay}, we expect that
the solution $X(t)$ is bounded in the norm $\Vert\cdot\Vert_{-\beta}$.
Then the integral (\re{FL}) converges and is analytic for $\Re\lambda>0$.
We will write $A$ and $v$ instead of $A_1$ and $v_1$ in all remaining part of the paper.
After the Fourier-Laplace transform  (\re{Avv}) becomes
\be\la{FLA}
  \lambda\tilde X(\lambda)=A\tilde X(\lambda)+X_0, \quad\Re\lambda>0.
\ee
Let us stress that (\re{FLA}) is equivalent to the Cauchy problem for the functions
$X(t)\in C_b([0,\infty);{\cal E}_{-\beta})$. Hence the solution $X(t)$ is given by
\be\la{FLAs}
  \tilde X(\lambda)=-(A-\lambda)^{-1}X_0,~~~~~~~~\Re\lambda>0
\ee
if the resolvent $R(\lambda)=(A-\lambda)^{-1}$ exists for $\Re\lambda>0$.

Let us comment on our following strategy in proving  the decay (\re{Zdec}). First, we will
construct the resolvent  $R(\lambda)$ for $\Re\lambda>0$ and prove that it is a continuous operator
in ${\cal E}_{-\beta}$. Then $\tilde X(\lambda)\in{\cal E}_{-\beta}$ and is an analytic function for
$\Re\lambda>0$. Second, we have to justify that there exist a (unique) function
$X(t)\in C([0,\infty);{\cal E}_{-\beta})$ satisfying (\re{FL}).

The analyticity of $\tilde X(\lambda)$ and Paley-Wiener arguments (see \ci{EKS}) should provide
the existence of a ${\cal E}_{-\beta}$ - valued distribution $X(t)$, $t\in\R$,
with a support in $[0,\infty)$. Formally,
\be\la{FLr}
  X(t)=\fr 1{2\pi}\int_\R e^{i\om t}\tilde X(i\om+0)d\om, ~~~~~~~~t\in\R.
\ee
However, to check the continuity of $X(t)$ for $t\ge 0$, we need additionally a bound for
$\tilde X(i\om+0)$ at large $|\om|$. Finally, for the time decay of $X(t)$, we need an additional
information on the smoothness and decay of $\tilde X(i\om+0)$. More precisely,
we should prove that the function $\tilde X(i\om+0)$ \\
i) is smooth outside $\om=0$ and $\om=\pm\mu$, where $\mu=\mu(v)>0$,
\\
ii) decays in a certain sense as $|\om|\to\infty$.
\\
iii)admits the Puiseux expansion at $\om=\pm\mu$.
\\
iv) is analytic at $\om=0$ if $X_0\in{\cal Z}_v:={\bf P}_v{\cal E}$ and $X_0\in{\cal E}_\beta$.
\\
Then the decay (\re{Zdec}) would follow from the Fourier-Laplace representation (\re{FLr}).

We will check the properties of type i)-iv)  only for the last two  components
$\tilde Q(\lambda)$ and $\tilde P(\lambda)$ of the vector $\tilde X(\lambda)
=(\tilde\Psi_1(\lambda),\tilde\Psi_2(\lambda),\tilde Q(\lambda),\tilde P(\lambda))$.
The properties provide the decay (\re{Zdec}) for the vector components
$Q(t)$ and $P(t)$ of the solution $X(t)$. Then for  the field components $\Psi_1(x,t)$ and
$\Psi_2(x,t)$ we will use wellknown properties of free  Schr\"odinger equation.



\setcounter{equation}{0}

\section{Constructing the Resolvent}

Here we construct the resolvent as a
bounded operator in ${\cal E}_{-\beta}$ for $\Re\lambda>0$.
We will write $(\Psi_1(y),\Psi_2(y), Q, P)$
instead of $(\tilde\Psi_1(y,\lambda),\tilde\Psi_2(y,\lambda),
\tilde Q(\lambda),\tilde P(\lambda))$
to simplify the notations. Then (\re{FLA}) reads
$$
(A-\lambda)\left(
\begin{array}{c}
\Psi_1 \\ \Psi_2 \\ Q \\ P
\end{array}
\right)=-\left(
\begin{array}{c}
\Psi_{01} \\ \Psi_{02} \\ Q_0 \\ P_0
\end{array}
\right).
$$
It is  the system of equations
\be\la{eq1}
\left.\begin{array}{r}
v\cdot\na\Psi_1(y)-(\Delta-m^2)\Psi_2(y)-Q\cdot\na\rho_2-\lambda\Psi_1(y)=-\Psi_{01}(y)
\\\\
(\Delta-m^2)\Psi_1(y)+v\cdot\na\Psi_2(y)+Q\cdot\na\rho_1-\lambda\Psi_2(y)=-\Psi_{02}(y)
\\\\
P-\lambda Q=-Q_0
\\\\
-\langle\na\Psi_j(y),\rho_j(y)\rangle+\langle\na\psi_{vj}(y),Q\cdot\na\rho_j(y)
\rangle-\lambda P=-P_0
\end{array}\right|~~~~~~~~~y\in\R^3.
\ee
{\it Step i)} Let us study the first two equations. In Fourier space they become
\be\la{F1}
  \left.\begin{array}{r}
  -(ikv+\lambda)\hat\Psi_1(k)+(k^2+m^2)\hat\Psi_2(k)=-\hat\Psi_{01}(k)-iQk\hat\rho_2,
  \\\\
  -(k^2+m^2)\hat\Psi_1(k)-(ikv+\lambda)\hat\Psi_2(k)=-\hat\Psi_{02}(k)+iQk\hat\rho_1.
  \end{array}\right|
\ee
Let us invert the matrix of the system and obtain
\begin{multline*}
  \!\!\!\!\!\left(\!\!
  \begin{array}{cc}
  -(ikv+\lambda) & k^2+m^2 \\
  -(k^2+m^2) & -(ikv+\lambda)
  \end{array}
  \right)^{-1}
  \!\!\!\!=[(ikv+\lambda)^2+(k^2+m^2)^2]^{-1}\!\!\left(\!\!
  \begin{array}{cc}
  -(ikv+\lambda) & -(k^2+m^2) \\
  k^2+m^2 & -(ikv+\lambda)
  \end{array}
  \right).
\end{multline*}
Taking the inverse Fourier transform we obtain the corresponding fundamental solution
\be\la{Green}
  G_{\lambda}(y)=\left(
  \begin{array}{cc}
  v\cdot\na-\lambda & \Delta-m^2 \\
  -\Delta+m^2       & v\cdot\na-\lambda
  \end{array}
  \right)g_{\lambda}(y),
\ee
where
\begin{multline}\la{dete}
 \!\!\! g_\lambda(y)=F^{-1}_{k\to y}\ds\fr{1}{(k^2+m^2)^2-(kv-i\lambda)^2}
  =F^{-1}_{k\to y}\ds\fr{1}{(k^2+m^2-kv+i\lambda)(k^2+m^2+kv-i\lambda)}
\end{multline}
Note that denominator in RHS \eqref{dete} does not vanish for $\Re\lambda>0,\,k\in\R^3$.
Moreover it  does not vanish for $\Re\lambda>0,\,k\in\C^3$  for sufficiently small $|\Im k|$.
Therefore $g_\lambda(y)$ decays exponentially by the Paley-Wiener arguments.
Let us compute the entries of matrix $G_{\lambda}$ explicitly:
\beqn\la{Gij}
  G^{11}_{\lambda}(y)&=&G^{22}_{\lambda}(y)
  =F^{-1}\ds\fr{-ikv-\lambda}{(k^2+m^2)^2-(kv-i\lambda)^2}\\
  \nonumber
  &=&F^{-1}_{k\to y}\Big(\ds\fr{1/2i}{k^2+m^2-kv+i\lambda}-\ds\fr{1/2i}{k^2+m^2+kv-i\lambda}\Big)
  =\ds\fr{e^{-\varkappa_+|y|-i\fr v2 y}}{8i\pi| y|}-\fr{e^{-\varkappa_-|y|+i\fr v2 y}}{8i\pi| y|},
  \\
  \nonumber
  G^{21}_{\lambda}(y)&=&-G^{12}_{\lambda}(y)
  =F^{-1}\ds\fr{k^2+m^2}{(k^2+m^2)^2-(kv-i\lambda)^2}\\
  \nonumber
  &=&F^{-1}_{k\to y}\Big(\ds\fr{1/2}{k^2+m^2-kv+i\lambda}+\ds\fr{1/2}{k^2+m^2+kv-i\lambda}\Big)
  =\ds\fr{e^{-\varkappa_+|y|-i\fr v2 y}}{8\pi| y|}+\fr{e^{-\varkappa_-|y|+i\fr v2 y}}{8\pi| y|},
\eeqn
where
\be\la{varkappa}
\varkappa_{\pm}=\sqrt{m^2-\fr{v^2}4\pm i\lambda},\quad\Re\varkappa_{\pm}>0.
\ee
This implies
\begin{lemma}\la{cres}
  i) The operator $G_{\lambda}$ with the integral kernel $G_{\lambda}(y-y')$,
  is continuous operator
  $H^1(\R^3)\oplus H^1(R^3)\to H^2(\R^3)\oplus H^2(R^3)$ for $\Re\lambda>0$.\\
  ii) The formulas \eqref{Gij} and  \eqref{varkappa} imply that for every fixed $y$,
  the matrix function $G_\lambda(y)$, $\Re\lambda>0$, admits an analytic continuation in $\lambda$
  to the Riemann surface of the algebraic function $\sqrt{\mu^2+\lambda^2}$
  with the branching points ~$\lambda=\pm i\mu$, where $\mu:=m^2-\fr{v^2}4$.
\end{lemma}
Thus, from (\re{F1}) and
 (\re{Green}) we obtain the convolution representation
\beqn\la{Psi}
  \Psi_1&=&-G^{11}_{\lambda}*\Psi_{01}-G^{12}_{\lambda}*\Psi_{02}
  -(G^{12}_{\lambda}*\na\rho_1)\cdot Q+(G^{11}_{\lambda}*\na\rho_2)\cdot Q,\\
  \nonumber
 \Psi_2&=&G^{12}_{\lambda}*\Psi_{01}-G^{11}_{\lambda}*\Psi_{02}
  -(G^{11}_{\lambda}*\na\rho_1)\cdot Q-(G^{12}_{\lambda}*\na\rho_2)\cdot Q.
\eeqn


\noindent{\it Step ii)}
Now we proceed to the last two equations (\re{eq1}):
\be\la{lte}
-\lambda Q+P=-Q_0,\quad\quad \langle\na\psi_{vj},
Q\cdot\na\rho_j\rangle-\langle\na\Psi_j,\rho_j\rangle-\lambda P=-P_0.
\ee
Let us rewrite equations (\re{Psi}) as
$\Psi_j=\Psi_{j}(Q)+\Psi_{j}(\Psi_{01},\Psi_{02})$, where
$$
  \Psi_1(\Psi_{01},\Psi_{02})=-G^{11}_{\lambda}*\Psi_{01}-G^{12}_{\lambda}*\Psi_{02},\quad
  \Psi_1(Q)=(-G^{12}_{\lambda}*\na\rho_1+G^{11}_{\lambda}*\na\rho_2)\cdot Q,
$$
$$
  \Psi_2(\Psi_{01},\Psi_{02})
  =G^{12}_{\lambda}*\Psi_{01}-G^{11}_{\lambda}*\Psi_{02},\quad\quad
  \Psi_2(Q)=-(G^{11}_{\lambda}*\na\rho_1+G^{12}_{\lambda}*\na\rho_2)\cdot Q.
$$
Then $\langle\na\Psi_j,\rho_j\rangle=\langle\na\Psi_{j}(Q),\rho_j\rangle+
\langle\na\Psi_{j}(\Psi_{01},\Psi_{02}),\rho_j\rangle$, and the last equation (\re{lte}) becomes
$$
\langle\na\psi_{vj},Q\cdot\na\rho_j\rangle-\langle\na\Psi_{j}(Q),
\rho_j\rangle-\lambda P=-P_0+\langle\na\Psi_{j}(\Psi_{01},\Psi_{02}),\rho_j\rangle=:-P_0'.
$$
First we compute the term
$$
  \langle\na\psi_{vj},Q\cdot\na\rho_j\rangle=\sum_{lj}\langle\na\psi_{vj},Q_l\pa_l\rho_j\rangle=
  \sum_{lj} \langle\na\psi_{vj},\pa_l\rho_j\rangle Q_l.
$$
Applying the Fourier transform $F_{y\to k}$, we have by the Parseval
identity and (\re{hpsiv}) that
\beqn\la{Lij}
  \sum_j\langle\pa_i\psi_{vj},\pa_l\rho_j\rangle&=&
  \sum_j\langle -ik_i\hat\psi_{vj},-ik_l\hat\rho_j\rangle\\
  \nonumber
  &=&\langle k_i\frac{-(k^2+m^2)\hat\rho_1+ikv\hat\rho_2}
  {(k^2+m^2)^2-(kv)^2},k_l\hat\rho_1\rangle
  +\langle k_i\frac{-ikv\hat\rho_1-(k^2+m^2)\hat\rho_2}
  {(k^2+m^2)^2-(kv)^2},k_l\hat\rho_2\rangle\\
  \nonumber\\
  \nonumber
  &=&-\int\fr{k_ik_l\Bigl((k^2+m^2)(|\hat\rho_1|^2+|\hat\rho_2|^2)+i(kv)
  (\hat\rho_1\overline{\hat\rho}_2-\hat\rho_2\overline{\hat\rho}_1)\Bigr)dk}
  {(k^2+m^2)^2-(kv)^2} =:-  L_{il}.
\eeqn
As the result, $\langle\na\psi_{vj},Q\cdot\na\rho_j\rangle=-LQ$, where $L$
is the $3\times3$ matrix with the matrix elements $L_{il}$.
Now let us compute the term $-\langle\na\Psi_j(Q),\rho_j\rangle=\langle\Psi_j(Q),\na\rho_j\rangle$.
One has
\begin{multline*}
  \langle\Psi_j(Q),\pa_i\rho_j\rangle
  =\sum\limits_l\Big(\langle -G^{12}_{\lambda}*\pa_l\rho_1
  + G^{11}_{\lambda}*\pa_l\rho_2,\pa_i\rho_1\rangle
  -\langle G^{11}_{\lambda}*\pa_l\rho_1+ G^{12}_{\lambda}*\pa_l\rho_2,\pa_i\rho_2\rangle\Big)Q_l\\
  =\sum\limits_l H_{il}(\lambda)Q_l,
\end{multline*}
and again by the Parseval identity we have
\beqn\la{Cij}
H_{il}(\lambda):&=&\langle -G^{12}_{\lambda}*\pa_l\rho_1+G^{11}_{\lambda}*\pa_l\rho_2,\pa_i\rho_1\rangle
  -\langle G^{11}_{\lambda}*\pa_l\rho_1+G^{12}_{\lambda}*\pa_l\rho_2,\pa_i\rho_2\rangle\\
  \nonumber\\
  \nonumber
  &=&\langle [(k^2+m^2)\hat\rho_1-(ikv+\lambda)\hat\rho_2]\hat g_{\lambda}k_l,k_i\hat\rho_1\rangle
  +\langle[(ikv+\lambda)\hat\rho_1+(k^2+m^2)\hat\rho_2]\hat g_{\lambda}k_l,k_i\hat\rho_2\rangle\\
  \nonumber\\
  \nonumber
  &=&\int\fr{k_ik_l\Bigl((k^2+m^2)(|\hat\rho_1|^2+|\hat\rho_2|^2)
  +(ikv+\lambda)(\hat\rho_1\overline{\hat\rho}_2-\hat\rho_2\overline{\hat\rho}_1)\Bigr)dk}{(k^2+m^2)^2-(kv-i\lambda)^2}
\eeqn
The matrix $H$ is well defined for $\Re\lambda>0$ since the denominator does not vanish. As the result,
$-\langle\na\Psi_j(Q),\rho_j\rangle=HQ$, where $H$ is the  matrix with matrix elements $H_{il}$.
Finally the equations (\re{lte}) become
\be\la{Mlam}
  {\cal M}(\lambda)\left(
  \begin{array}{c}
  Q \\ P
  \end{array}
  \right)=\left(
  \begin{array}{c}
  Q_0 \\ P_0'
  \end{array}
  \right),\,\,{\rm where}\,\,{\cal M}(\lambda)=\left(
  \begin{array}{cc}
  \lambda E  & -E \\
  L-H(\lambda) & \lambda E
  \end{array}
  \right),
\ee
\\
 Assume for a moment that the matrix ${\cal M}(\lambda)$ is invertible (later we will prove this). Then we obtain
\be\la{QP1}
  \left(
  \begin{array}{c}
  Q \\ P
  \end{array}
  \right)={\cal M}^{-1}(\lambda)\left(
  \begin{array}{c}
  Q_0 \\ P_0'
  \end{array}
  \right),~~~~~~~~~\Re\lambda>0.
\ee
Finally, formula (\re{QP1}) and formulas (\re{Psi}), where $Q$ is expressed from (\re{QP1}),
give the expression of the resolvent $R(\lambda)=(A-\lambda)^{-1}$, $\Re\lambda>0$\\
\begin{lemma}\la{cmf}
  The matrix function ${\cal M}(\lambda)$ (respectively, ${\cal M}^{-1}(\lambda)$),
  $\Re\lambda>0$ admits
  an analytic (respectively meromorphic) continuation to the Riemann surface of the function\\
  $\sqrt{\mu^2+\lambda^2}$, $\lambda\in\C$.
\end{lemma}
\begin{proof}
  The analytic continuation of  ${\cal M}(\lambda)$,  exists by Lemma \re{cres} ii)
  and the
  convolution expressions in (\re{Cij}) by  (\re{ro}). The inverse matrix is then
  meromorphic since it exists for large $\Re\lambda$. The latter follows from (\re{Mlam})
  since $H(\lambda)\to 0$, $\Re\lambda\to\infty$, by  (\re{Cij}).
\end{proof}
\setcounter{equation}{0}
\section{Analyticity in the Half-Plane}
Here we prove the following
\begin{prop}\la{analyt}
  The operator-valued function $R(\lambda):{\cal E}\to{\cal E}$ is analytic for $\Re\lambda>0$.
\end{prop}
\begin{proof}
  It is sufficient to prove that the operator $A-\lambda:{\cal E}\to {\cal E}$ has a bounded inverse
  operator for $\Re\lambda>0$. Let us recall, that  $A=A_{v,v}$ where $|v|<2m$.\medskip \\
  {\it Step i)} Let us prove that Ker$\,(A-\lambda)=0$ for $\Re\lambda>0$. Indeed, let us assume
  that $X_{\lambda}=(\Psi_{\lambda 1},\Psi_{\lambda 2},Q_{\lambda},P_{\lambda})\in{\cal E}$
  satisfies $(A-\lambda)X_{\lambda}=0$, that is $X_{\lambda}$ is a solution to (\re{eq1}) with
  $\Psi_{01}=\Psi_{02}=0$ and $Q_0=P_0=0$. We have to prove that $X_{\lambda}=0$.

  First let us check that $P_{\lambda}=0$. Indeed, the trajectory
  $X:=X_{\lambda}e^{\lambda t}\in C(\R,{\cal E})$
  is the solution to the equation $\dot X=AX$ that is (\re{line}) with $w=v$.
  Then ${\cal H}_{v,v}(X(t))$ grows exponentially by  (\re{H0vv}). This growth contradicts to
  the conservation of ${\cal H}_{v,v}$, which follows from  Lemma \re{haml} ii) since
  $X(t)\in C^1(\R,{\cal E}^+)$. The latter inclusion follows from Lemma \re{cres} since
  $(\Psi_{\lambda 1},\Psi_{\lambda 2})$ satisfies the equations (\re{Psi})
  with $\Psi_{01}=\Psi_{02}=0$ and $Q=Q_\lambda$.

  Now $\lambda Q_\lambda=P_\lambda=0$ by the third equation of (\re{eq1}), hence $Q_{\lambda}=0$ since
  $\lambda\ne 0$. Finally, $\Psi_{\lambda 1}=0$, $\Psi_{\lambda 2}=0$ by the equations (\re{Psi})
  with $Q=Q_\lambda=0$.
  \medskip \\
  {\it Step ii)}
 Let us represent $A-\lambda =A_0+T$, where
\[  A_0\!=\!\left(
   \begin{array}{cccc}
   v\cdot\na-\lambda   & -(\Delta-m^2)  & 0 & 0 \\
   \Delta-m^2 & v\cdot\na-\lambda & 0 & 0 \\
    0 &  0 & \!\!\!\!-\lambda  & 0 \\
    0 &  0 & 0 &\!\!\!\! -\lambda
   \end{array}
   \right),    \;
   T\!=\!\left(
   \begin{array}{cccc}
   0 & 0 & \!\!\!\!\!\!\!-\cdot\na\rho_2  & 0 \\
   0 & 0 &            \cdot\na\rho_1& 0 \\
   0 & 0 & 0 & E \\
   \langle\cdot,\na\rho_1\rangle & \langle\cdot,\na\rho_1\rangle &
   \langle\na\psi_{vj},\cdot\na\rho_j\rangle & 0
   \end{array}
   \right).    \]
The operator $T$ is finite-dimensional, and the operator $A_0^{-1}$ is bounded
in ${\cal E}$ by Lemma \re{cres}. Finally, $A-\lambda =A_0(I+A_0^{-1}T)$, where
$A_0^{-1}T$ is a compact operator. Since we know that Ker$\,(I+A_0^{-1}T)=0$,
the operator $(I+A_0^{-1}T)$ is invertible by Fredholm theory.
\end{proof}

\begin{cor}
The matrix ${\cal M}(\lambda)$ of (\re{Mlam}) is invertible for $\Re\lambda>0$.
\end{cor}
\setcounter{equation}{0}

\section{Regularity on the Imaginary Axis}

First, let us describe  the continuous spectrum of the operator $A=A_{v,v}$
on the imaginary axis. By definition, the continuous spectrum corresponds to
$\om\in\R$, such that the resolvent $R(i\om+0)$ is not a bounded operator
in ${\cal E}$. By the formulas (\re{Psi}), this is the case when the Green function
$G_\lambda(y-y')$ loses the exponential decay.
Thus, $i\om$ belongs to the continuous spectrum if
\[   |\om|\ge\mu=m^2-v^2/4.   \]
By Lemma \re{cmf}, the limit matrix
\be\la{M}
 {\cal M}(i\om):={\cal M}(i\om+0)=\left(
  \begin{array}{cc}
  i\om E        & -E      \\
  L  -H(i\om+0) & i\om E
  \end{array}
  \right), ~~~~~~~~~~\om\in\R,
\ee
exists, and its entries are continuous functions of $\om\in\R$, smooth for $|\om|<\mu$
and $|\om|>\mu$. Recall that the point $\lambda=0$ belongs to the discrete spectrum of
the operator $A$ by Lemma \re{ceig} i), hence ${\cal M}(i\om+0)$
(probably)  also is not invertible at $\om=0$.
\begin{prop}\la{regi}
  Let $\rho$ satisfy the condition\eqref{ro} and the Wiener condition (\re{W}), and $|v|<2m$.
  Then the limit matrix ${\cal M}(i\om+0)$ is invertible for $\om\ne 0$, $\om\in\R$.
\end{prop}
{\bf Proof } We will consider  separately three cases
$0<|\om|<\mu$, $\om= \mu$, and $|\om|>\mu$. We can assume that
$v=(|v|,0,0)$.
Let us denote $F(\om):=-L+H(i\om+0)$,
 $M=m^2+k^2$, $a=|\hat\rho_1|^2+|\hat\rho_2|^2$,
$b=i(\hat\rho_1\overline{\hat\rho}_2-\hat\rho_2\overline{\hat\rho}_1)$.
Then the entries of the
matrix $F$ become
\beqn\la{decomp}
  F_{ij}&=&\int k_ik_j~dk \Bigg[Ma\left(\frac 1{M^2-(|v|k_1+\omega)^2}
  -\frac 1{M^2-(|v|k_1)^2}\right)\\
  \nonumber
  &+& b\left(\frac{|v|k_1+\omega}
  {M^2-(|v|k_1+\omega)^2}-\frac{|v|k_1}{M^2-(|v|k_1)^2}\right)\Bigg]\\
  \nonumber
  &=&\int\fr{k_ik_jdk}2 \Bigg[a\left(\frac 1{M-|v|k_1-\omega}+\frac 1{M+|v|k_1+\omega}
  -\frac 1{M-|v|k_1}-\frac 1{M+|v|k_1}\right)\\
  \nonumber
  &+& b\left(\frac 1{M-|v|k_1-\omega}-\frac 1{M+|v|k_1+\omega}
  -\frac 1{M-|v|k_1}+\frac 1{M+|v|k_1}\right)\Bigg].
\eeqn
Since $a$ is even, and $b$ is  odd we obtain that
\be\la{decomp1}
  F_{ij}=\fr 12\int~dk_2dk_3\int\limits_0^{+\infty} k_ik_j~dk_1\Big[af_1+bf_2\Big]
\ee
where
\beqn\label{cM}
  f_1:&=&\fr1{M-|v|k_1-\om}+\fr1{M+|v|k_1+\om}+\fr1{M+|v|k_1-\om}+\fr1{M-|v|k_1+\om}\\
  \nonumber
  &-&\fr2{M-|v|k_1}-\fr2{M+|v|k_1},\\
  \nonumber
  f_2:&=&\fr1{M-|v|k_1-\om}-\fr1{M+|v|k_1+\om}+\fr1{M-|v|k_1+\om}-\fr1{M+|v|k_1-\om}\\
  \nonumber
  &-&\fr2{M-|v|k_1}+\fr2{M+|v|k_1},
\eeqn

Then by (\re{M})
\begin{multline}\la{detM}
{\rm det}\,{\cal M}(i\om)
={\rm det}\left(
\begin{array}{cccccc}
   i\om     &     0      &     0     &    -1     &    0     &   0    \\
    0       &    i\om    &     0     &     0     &   -1     &   0    \\
    0       &     0      &   i\om    &     0     &    0     &  -1    \\
 -F_{11}    &  -F_{12}   &  -F_{13}  &    i\om   &    0     &   0    \\
 -F_{12}    &  -F_{22}   &  -F_{23}  &     0     &   i\om   &   0    \\
 -F_{13}    &  -F_{23}   &  -F_{33}  &     0     &    0     &  i\om
\end{array}
\right)\\
=-\omega^6-\omega^4\sum\limits_{j=1}^3F_{jj}
-\omega^2\sum\limits_{i<j}(F_{ii}F_{jj}-F_{ij}^2)
-{\rm det}\left(
\begin{array}{ccc}
 F_{11}    &  F_{12}   &  F_{13} \\
 F_{12}    &  F_{22}   &  F_{23} \\
 F_{13}    &  F_{23}   &  F_{33}
 \end{array}
\right)
\end{multline}
since $F_{ij}=F_{ji}$.

\noindent

{\bf I.} First, let us consider the case $0<|\om|<\mu$. Then the invertibility
of ${\cal M}(i\om)$ follows from
\begin{lemma}\la{lnW}
  For $0<|\om|<\mu$,
  the matrix $F$ is positive definite.
\end{lemma}
\begin{proof}
  First, let us note that all denominators in \eqref{cM} are positive for
  $|\om|<\mu=m^2-v^2/4,\; |v|<2m$. Indeed,
  $$(m^2+k^2)^2-(\om+|v|k_1)^2=((k-v/2)^2+m^2-\frac{v^2}4-\om)((k+v/2)^2+m^2-\frac{v^2}4+\om)>0$$
  Second, $f_1>f_2\ge 0$ if $|v|<2m$ and $0<|\om|\le\mu$. This is proved in Appendix C.\\
  Finally, the Wiener condition implies
  \be\la{pd}
  a\pm b=|\hat\rho_1(k)\mp i\hat\rho_2(k)|^2>0,\quad\forall k\in\R^3.
  \ee
  Therefore $af_1+bf_2>0$ and \eqref{decomp1}
   is the integral of the symmetric nonnegative definite
  matrix $k\otimes k=(k_ik_j)$ with a positive weight.
  Hence, the matrix $F$ is positive definite.
 \end{proof}

\medskip

\noindent {\bf II.} $\om=\pm\mu$. Let us consider for example
$\om=\mu=m^2-\fr{|v|^2}4$. In this case (\re{Cij}) reads:
$$
H_{ij}(i\mu)=\int\fr{k_ik_j(Ma-(kv+\mu)b)dk}
{\Big((k_1-\frac{|v|}2)^2+k_2^2+k_3^2\Big)\Big((k_1+\frac{|v|}2)^2+k_2^2+k_3^2+2\mu\Big)}.
$$
Now the integrand has a unique singular point. The singularity is
integrable, hence ${\rm det}\,{\cal M}(i\om)$ also is negative by the
representations (\re{detM}). Hence, the matrix ${\cal M}(i\mu)$ is also invertible.
\smallskip

\noindent {\bf III.} $|\om|>\mu$. Here we apply an other arguments.
Now the invertibility  of ${\cal M}(i\om)$ follows from  (\re{detM}) by
the following lemma (cf. \cite{ikv})
\begin{lemma}\la{lW}
  If (\re{W}) holds and $\om>\mu$ ($\om<-\mu$), then the matrix $\Im F(\om)$
  is negative (positive) definite.
\end{lemma}
\begin{proof}
  We consider the case $\om>\mu$ (the case $\om<-\mu$ can be treated similarly).
  Let us calculate the imaginary part of $F_{ij}$. Since
  $F_{ij}=H_{ij}(i\om+0)-L_{ij}$ and $L_{ij}$ is real, we will
  consider only $H_{ij}(i\om+0)$. For $\ve>0$ we have
  \beqn\la{Hjjlim}
    H_{ij}(i\om+\ve)
    &=&\int\fr{k_ik_j(Ma+(kv+\om-i\ve)b)dk}{M^2-(kv+\om-i\ve)^2}
    =\fr 12\int\fr{k_ik_j(a+b)dk}{M-kv-\om+i\ve}\\
    \nonumber
    &+&\fr 12\int\fr{k_ik_j(a-b)dk}{M+kv+\om-i\ve}
    = H^1_{ij}(i\om+\ve)+ H^2_{ij}(i\om+\ve).
  \eeqn
  It  suffices to consider only the first summand in \eqref{Hjjlim}, since the second
  summand is real for $\ve=0$.
  Consider the denominator
  $$ \hat D_{\ve}(k)=k^2+m^2-kv-\om+i\ve.$$
  Note that $\hat D_{0}(k)=0$ on an ellipsoid $T_{\om}$, where
  $$ T_{\om}=\{k:|k-\fr{v}2|=R:=\sqrt{\om-\mu}\},$$
  Then  the Plemelj formula for $C^1$-functions implies that
  \be\la{ImHjj}
    \Im H^1_{ij}(i\om+0)=-\fr{\pi}2\int_{T_{\om}}\fr{k_ik_j(a+b)}
    {|\na\hat D_{0}(k)|}dS,
  \ee
  where $dS$ is the element of the surface area.
  Hence, the matrix $\Im H^1(i\om+0)$ is negative definite by (\re{pd}).
\end{proof}
Now let us prove that the limit matrix ${\cal M}(i\om +0)$ is invertible.  Recall that
\begin{equation*}
 {\cal M}(i\om +0)=
\left(
 \begin{array}{ll}
    i\om E            &  -E   \\
  - F(i\om +0)         &   i\om E
 \end{array}
\right)
\end{equation*}
Then the equation
\begin{equation*}
{\cal M}(i\om +0)\left(
\begin{array}{c}
Q \\ P
\end{array}
\right)=0
\end{equation*}
becomes
\be\la{MPQs}
i\om Q -P=0,\quad -FQ+i\om P=0.
\ee
Then $(F+\om^2)Q=0$, which implies $Q=0$  and then $P=0$ since the matrix $\Im F$
is  negative  definite for $\om>\mu$.
This completes the proofs of  the Proposition 15.1.
\begin{cor}\la{creg}
  Proposition \re{regi} implies that the matrix ${\cal M}^{-1}(i\om)$ is smooth in $\om\in\R$
  outside three points $\om=0,\pm \mu$.
\end{cor}

\setcounter{equation}{0}
\section{Singular Spectral Points}

Let us recall that the formula (\re{QP1}) expresses the Fourier-Laplace transforms
$\tilde Q(\lambda),\tilde P(\lambda)$. Hence, the components are given by the Fourier integral
\be\la{QP1i}
  \left(
  \begin{array}{c}
  Q(t) \\ P(t)
  \end{array}
  \right)=\ds\fr 1{2\pi}\int e^{i\om t}{\cal M}^{-1}(i\om+0)\left(
  \begin{array}{c}
  Q_0 \\ P_0'
  \end{array}
  \right)d\om
\ee
if it converges in the sense of distributions. The Corollary \re{creg} alone is not sufficient for the proof
of the convergence and decay of the integral. Namely, we need an additional information about a regularity
of the matrix ${\cal M}^{-1}(i\om)$ at its singular points  $\om=0,\pm \mu$, and some bounds at $|\om|\to\infty$.
We will analyze all the points separately.
\medskip\\
{\bf I.}
First we consider the points $\pm \mu$.
\begin{lemma} \la{Pui}
  The matrix  ${\cal M}^{-1}(i\om)$ admits the following Puiseux expansion
  in a neighborhood of $\pm \mu$: there exists an $\ve_\pm>0$ s.t.
  \be\la{mom}
    {\cal M}^{-1}(i\om)=\sum_{k=0}^{\infty}R_k^\pm(\om\mp\mu)^{k/2},\quad |\om\mp\mu|<\ve_\pm,\quad\om\in\R.
  \ee
\end{lemma}
\begin{proof}
  It suffices to prove a similar expansion  for ${\cal M}(i\om)$. Then  (\re{mom}) holds also
  for ${\cal M}^{-1}(i\om)$,
  since the matrices ${\cal M}(\pm i\mu)$ are invertible. The asymptotics for ${\cal M}(i\om)$ holds by
  the convolution representation  (\re{Cij})
 \begin{multline}\la{Cjjj}
  H_{ij}(\lambda)=-\langle G^{12}_{\lambda}*\pa_l\rho_1+G^{11}_{\lambda}*\pa_l\rho_2,\pa_i\rho_1\rangle
  -\langle  G^{11}_{\lambda}*\pa_l\rho_1+G^{12}_{\lambda}*\pa_l\rho_2,\pa_i\rho_2\rangle.
 \end{multline}
since  $G^{ij}_{\lambda}$ admit the corresponding Puiseux expansions by the formula (\re{Gij}).
\end{proof}
{\bf II.}
Second, we study the asymptotic behavior of ${\cal M}^{-1}(\lambda )$ at
infinity. Let us recall that ${\cal M}^{-1}(\lambda )$ was originally defined for
$\Re\lambda >0,$ and admits a meromorphic continuation to the Riemann
surface of the function $\sqrt{m^2-\fr{v^2}4+ i\lambda}$ (see Lemma \ref{cmf}).

\begin{lemma}\la{162}
There exist a matrix $R_{0}$ and a matrix-function $R_{1}(\om )$,  such that
$$
{\cal M}^{-1}(i\om)=\fr {R_0}\om +R_1(\om ),~~~|\om|\ge\mu+1,~~~~~~~~~\om \in \R,
$$
where, for every $k=0,1,2,..$
\be\label{min}
|\pa_\om^k R_1(\om )|\leq \fr{C_k}{|\om |^2},
~~~~~~~~~~~~|\om|\ge\mu+1,~~~~~~~~~\om \in \R.
\ee
\end{lemma}
\begin{proof}
The structure (\ref{M}) of the matrix ${\cal M}(i\om )$ provides that it suffices to prove the
following estimate for the elements of the matrix $H(i\om ):=H(i\om+0)$:
\be \la{minC}
  |\pa _\om^k H_{jj}(i\om )|\leq \fr{C_k}
  {|\om|}, ~~~~~~~~~\om \in \R,~~~|\om|\ge\mu+1,~~~~j=1,2,3.
\ee
Note, that
$$ G^{11}_{\lambda}*f=\ds\fr 1{2i}(D^{-1}_1(\lambda)f-D^{-1}_2(\lambda)f),\quad
 G^{12}_{\lambda}*f=\ds\fr 12(D^{-1}_1(\lambda)f+D^{-1}_2(\lambda)f),$$
where
$$ D_1(\lambda)=-\Delta+m^2-iv\cdot\nabla+i\lambda,\quad
   D_2(\lambda)=-\Delta+m^2+iv\cdot\nabla-i\lambda,\;\Re\lambda>0,$$
and $D_j^{-1}(\lambda),\;j=1,2$ are   bounded operators in $L^2(\R^3)$.
The estimate (\ref{minC}) immediately follows from a more general bound
\be \la{minCg}
\Vert\pa _\om^k D_j^{-1}(i\om+0)f\Vert_{L^2_{-\sigma}}\le \fr {C_k(R)}{|\om|}
\Vert f\Vert_{L^2_\sigma},~~~~~~~~~\om \in \R,~~~|\om|\ge \mu+1
\ee
which holds for  $\sigma>3/2$. Namely, (\ref{minC}) follows by (\re{ro}) from (\ref{minCg})
applied to the functions $f(y)=\pa_l\rho_j(y)\in  L^2_\sigma$.\\
The bound (\ref{minCg}) is proved in \ci[the bound (A.2')]{Ag} (see also \ci[Thm 8.1]{JK}).
\end{proof}

{\bf III.}
Finally, we consider the point  $\om=0$ which is most singular. The point is an isolated pole
of a finite degree by Lemma \re{cmf}, hence the Laurent expansion holds,
\be\la{Lor}
  {\cal M}^{-1}(i\om)=\sum_{k=0}^n  M_k\om^{-k-1}+{\cal H}(\om),~~~~~~~|\om|<\ve_0,
\ee
where $M_k$ are $6\times 6$ complex matrices, $\ve_0>0$, and ${\cal H}(\om)$ is
an analytic matrix-valued function for complex $\om$ with $|\om|<\ve_0$.


\section{Time Decay of the Vector Components}

Here we prove the decay (\re{Zdec}) for the components $Q(t)$ and $P(t)$.
\begin{lemma}\la{171} (cf. \cite{ikv})
Let $X_0\in {\cal Z}_{v,\beta}$. Then $Q(t)$, $P(t)$ are continuous and
\be\la{decQP}
|Q(t)|+|P(t)|\le \ds\fr {C(\rho, \overline v,d_0)}{(1+|t|)^{3/2}},
~~~~~~~t\ge 0.
\ee
\end{lemma}
\begin{proof}
The expansions (\re{mom}),
(\re{min})  and (\re{Lor}) imply
the convergence of the Fourier integral (\re{QP1i})
in the sense of distributions to a continuous function
of $t\ge 0$. Let us prove (\re{decQP}).
First let us note that the condition  $X_0\in {\cal Z}_{v,\beta}$ implies that
the whole trajectory $X(t)$ lies in  ${\cal Z}_{v,\beta}$. This follows
from the invariance of the space ${\cal Z}_{v,\beta}$ under the generator
$A_{v,v}$ (cf. Remark \re{rZ}).
Note that for $X_0$ not belonging to ${\cal Z}_{v,\beta}$
the components  $Q(t)$ and $P(t)$ may contain non-decaying terms
which
correspond to the singular point $\om=0$.
Indeed, we
know that the linearized dynamics admits the secular solutions without decay,
see (\re{secs}).
The formulas (\re{inb}) give
the corresponding components $Q_S(t)$ and $P_S(t)$ of the
secular solutions,
\be\la{secQPs}
  \left(\begin{array}{c}Q_S(t)\\P_S(t)\end{array}\right)
  =\sum_1^3 C_j
  \left(\begin{array}{c}e_j\\0\end{array}\right)
  +\sum_1^3 D_j\Bigg[\left(\begin{array}{c}e_j\\0\end{array}\right)t
  +\left(\begin{array}{c}0\\ e_j\end{array}\right)\Bigg].
\ee
We will show that  the symplectic orthogonality condition leads to (\re{decQP}).
Let us split the Fourier integral (\re{QP1i}) into three terms using the partition of unity
$\zeta_1(\om)+\zeta_2(\om)+\zeta_3(\om)=1$, $\om\in\R$:
\beqn\la{QP1i3}
\left(
  \begin{array}{c}
  Q(t) \\ P(t)
  \end{array}
  \right)&=&\ds\fr 1{2\pi}\int e^{i\om t}(\zeta_1(\om)+\zeta_2(\om)+\zeta_3(\om))
  {\cal M}^{-1}(i\om+0)\left(
  \begin{array}{c}
   Q_0 \\ P'_0
  \end{array}
  \right)d\om
  \nonumber\\
  \nonumber\\
  &=&I_1(t)+I_2(t)+I_3(t),
\eeqn
where the functions $\zeta_k(\om)\in C^\infty(\R)$ are supported by
\be\la{zsup}
\left.\begin{array}{rcl}
\supp \zeta_1&\subset& \{\om\in\R:\ve_0/2<|\om|<\mu+2\}
\\
\\
\supp \zeta_2&\subset& \{\om\in\R:|\om|>\mu+1\}
\\
\\
\supp \zeta_3&\subset& \{\om\in\R:|\om|<\ve_0\}
\end{array}\right|
\ee
Then
\smallskip\\
i) The function $I_1(t)$ decays, like $(1+|t|)^{-3/2}$, by the Puiseux expansion (\re{mom}).
\\
ii) The function $I_2(t)$ decays faster than any power of $t$ due to Proposition \re{162}.
\\
iii) Finally, the function $I_3(t)$ generally does not decay if $n\ge 0$ in the Laurent expansion
(\re{Lor}). Namely, the contribution of the analytic function ${\cal H}(\om)$ decays
faster than any power of $t$. On the other hand, the contribution of the Lorent series,
\be\la{Fin}
  \left(
  \begin{array}{c}
  Q_L(t) \\ P_L(t)
  \end{array}
  \right):=\ds\fr 1{2\pi} \int e^{i\om t}\zeta_3(\om)\sum_{k=0}^n M_k(\om-i0)^{-k-1}
  \left(
  \begin{array}{c}
  Q_0 \\ P_0'
  \end{array}
  \right)d\om,~~~~~~~t\in\R,
\ee
is a polynomial function of $t\in\R$, of a degree $\le n$, modulo functions decaying faster
than any power of $t$. Let us note that the formula (\re{secQPs}) gives an example
of the polynomial functions appeared from (\re{Fin}).

We have to show that the symplectic orthogonality condition eliminates the polynomial functions.
Our main difficulty is that we do not know anything about the order $n$ of the pole
and the Lorent coefficients $M_k$ of the matrix ${\cal M}^{-1}(i\om)$ at $\om=0$.

Our crucial observation is the following:
\smallskip\\
a) The components  (\re{secQPs}), of the secular solutions, form a linear space ${\cal L}_S$
of the dimension dim~${\cal L}_S= 6$.
\\
b) The polynomial functions  from (\re{Fin}) belong to a linear space ${\cal L}_L$
of the dimension dim~${\cal L}_L\le 6$, since  $(Q_0,P'_0)\in\R^6$.
\\
c)
${\cal L}_S\subset{\cal L}_L$ since all the functions (\re{secQPs}) admits the representation (\re{Fin}).
The latter follows from the fact that the secular solutions  (\re{secs}) can be reproduced by
our calculations with the Laplace transform.
\smallskip\\
Therefore, we conclude that
\be\la{LLL}
{\cal L}_L ={\cal L}_S.
\ee
It remains to note that the secular solutions are forbidden since $X_0\in {\cal Z}_{v,\beta}$.
Hence, the polynomial terms in (\re{Fin}) vanish that implies the decay (\re{decQP}).

More precisely, we know that$X(t)={\bf P}_vX(t)$ for all $t\in\R$. On the other hand,
the identity (\re{LLL}) implies that $X(t)$ can be corrected by a secular
solution $X_S(t)$ s.t. the corresponding components $Q_\Delta(t)$ and $P_\Delta(t)$,
of the difference $\Delta(t):=X(t)-X_S(t)$, decay. Hence, the components  $Q(t)$ and $P(t)$,
of $X(t)={\bf P}_vX(t)={\bf P}_v[X(t)-X_S(t)]$, also decay.
\end{proof}

\setcounter{equation}{0}



\setcounter{equation}{0}

\section{Time Decay of Fields}
Here we prove the decay  of the field components $\Psi_1(x,t),\Psi_2(x,t)$ corresponding
to (\re{Zdec}). The first  two equations of (\re{Avv}) may be written as one equation:
\be\la{linf}
  i\dot \Psi(t)=(-\Delta+m^2+iv\cdot\nabla)\Psi-Q(t)\cdot\na\rho,
\ee
where $\Psi(t)=\Psi_1(\cdot,t)+i\Psi_2(\cdot,t)).$
By Lemma \re{171} we know that $Q$ is  continuous function of $t\ge 0$, and
\be\la{linb}
  |Q(t)|\le \ds\fr {C(\rho, \overline v,d_0)}{(1+|t|)^{3/2}},~~~~t\ge 0.
\ee
Hence, the Proposition \re{lindecay} is reduced now to the following
\begin{prop} \la{pfi}
i) Let a function $Q(t)\in C([0,\infty);\R^3)$,
and $\Psi_0\in H^1_\beta$. Then the equation (\re{linf}) admits a unique
solution $\Psi(t)\in C([0,\infty); H^1_\beta)$ with the initial condition
$\Psi(0)=\Psi_0$.
\\
ii) If $\Psi_0\in H^1_{\beta}$ and the decay (\re{linb}) holds,
the corresponding fields also decay uniformly in $v$:
\be\la{lins}
\Vert \Psi(t)\Vert_{1,-\beta}
\le \ds\fr{C(\rho,\overline v,d_0,\Vert\Psi_0\Vert_{1,\beta})}{(1+|t|)^{3/2}},~~~~t\ge 0,
\ee
for $|v|\le \overline v$ with any $\overline v\in (0,2m)$.
\end{prop}
\begin{proof}
 The statements follow from the Duhamel representation
\be\la{Duh}
\Psi(t)=W(t)\Psi_0-\int_0^tW(t-s) Q(s)\cdot\na\rho~ds,
~~~~~~t\ge 0, \ee
where $W(t)$ is the dynamical group (propagator) of the free equation
$$
 i\dot\Psi(t)=(-\Delta+m^2+iv\cdot\nabla)\Psi(t).
$$
\begin{lemma}
Let $|v|\le \overline v$ with any $\overline v\in (0,2m)$. Then
for $\Psi_0\in H^1_\beta$.
\be\la{tS}
\Vert W(t)\Psi_0\Vert_{1,-\beta}\le C(\ov
v)(1+|t|)^{-3/2}\Vert\Psi_0\Vert_{1,\beta},\;~~~~t\ge 0.
\ee
\end{lemma}
\begin{proof}
Note that $W(t)\Psi_0= e^{-i(m^2-{|v|^2}/4)t}e^{i{v}x/2}\Phi(t), $
where $\Phi(t)$ is a solution to free Schr\"o\-din\-ger Equation
$$i\dot \Phi(t)=-\Delta\Phi(t),\quad \Phi(0)=e^{i{v}x/2}\Psi_0.$$
It is wellknown  that  $\Phi(t)$ satisfies the estimate
$\Vert\Phi(t)\Vert_{1,-\beta}\le C(1+|t|)^{-3/2}\Vert\Phi(0)\Vert_{1,\beta}$,
$t\ge 0$ (see for example \cite{JK}).
\end{proof}

 Now \eqref {lins}  follows from  the condition
(\re{linb}), and the Duhamel representation (\re{Duh}).
\end{proof}
\setcounter{equation}{0}
\section{Appendix}

\subsection*{A. Solitary waves}

Let us to check the last equation of \eqref{stfch}:
\be\la{last}
0=\int\bigl(\nabla\psi_{v1}(y)\rho_1(y)+\nabla\psi_{v2}(y)\rho_2(y)\bigr)dy.
\ee
Let us transfer to the Fourier representation. Set
$$
   \hat\psi(k):=(2\pi)^{-3/2}\int e^{ikx}\psi(x)dx.
$$
It is easy to compute that
\be\la{system}
   -ikv\hat\psi_{v1}+(k^2+m^2)\hat\psi_{v2}=-\hat\rho_2,\quad
   (k^2+m^2)\hat\psi_{v1}+ikv\hat\psi_{v2}=-\hat\rho_1.
\ee
Therefore
\be\la{hpsiv}
   \hat\psi_{v1}(k)=\fr{-(k^2+m^2)\hat\rho_1(k)
   +ikv\hat\rho_2(k)}{(k^2+m^2)^2-(kv)^2},\quad
   \hat\psi_{v2}(k)=\fr{-ikv\hat\rho_1(k)-(k^2+m^2)
   \hat\rho_2(k)}{(k^2+m^2)^2-(kv)^2}.
\ee
By Parseval identity \eqref{last} becames
\begin{multline*}
0=\int k_j\bigl(\hat\psi_{v1}\overline{\hat\rho}_1
+\hat\psi_{v2}\overline{\hat\rho}_2\bigr)dk
=\int\frac{k_j\bigl[-(k^2+m^2)(|\hat\rho_1|^2+|\hat\rho_2|^2)
+ikv(\hat\rho_2\overline{\hat\rho}_1-\hat\rho_1\overline{\hat\rho}_2)\bigr]dk}
{(k^2+m^2)^2-(kv)^2},
\end{multline*}
which is true, since the integrand is odd.

\subsection*{B. Computing $\Om(\tau_i,\tau_j)$}
Let us to justify the formulas  (\re{Omega})-(\re{alpha}) for the matrix $\Om$.
For $j,l=1,2,3$ one has from (\re{inb}) and (\re{OmJ})
\be\la{jl}
   \Om(\tau_j,\tau_l)=\langle\pa_j\psi_{v1},\pa_l\psi_{v2}\rangle-
   \langle\pa_j\psi_{v2},\pa_l\psi_{v1}\rangle,
\ee
\be\la{jp3lp3}
   \Om(\tau_{j+3},\tau_{l+3})=\langle\pa_{v_j}\psi_{v1},\pa_{v_l}\psi_{v2}\rangle-
   \langle\pa_{v_j}\psi_{v2},\pa_{v_l}\psi_{v1}\rangle,
\ee
and
\be\la{Wx}
   \Om(\tau_{j},\tau_{l+3})
   =-\langle\pa_{j}\psi_{v1},\pa_{v_l}\psi_{v2}\rangle+
   \langle\pa_{j}\psi_{v2},\pa_{v_l}\psi_{v1}\rangle+e_j\cdot e_l.
\ee
Differentiating \eqref{system} we get
\be\la{derv}
   \pa_{v_j}\hat\psi_{v1}=\fr{k_jkv\hat\psi_{v1}
   -ik_j(k^2+m^2)\hat\psi_{v2}}{(k^2+m^2)^2-(kv)^2},\quad
   \pa_{v_j}\hat\psi_{v2}=\fr{ik_j(k^2+m^2)\hat\psi_{v1}
   +k_jkv\hat\psi_{v2}}{(k^2+m^2)^2-(kv)^2},\;j=1,2,3.
\ee
Then for $j,l=1,2,3$ we obtain from (\re{jl}) by the  Parseval identity that
\be\la{jl123}
   \Om(\tau_j,\tau_l)=\int k_jk_l~dk(\hat\psi_{v1}\overline{\hat\psi}_{v2}
   -\hat\psi_{v2}\overline{\hat\psi}_{v1})= 0,
\ee
since the function $\hat\psi_{vc}=
\hat\psi_{v1}\overline{\hat\psi}_{v2} -\hat\psi_{v2}\overline{\hat\psi}_{v1}$ is odd.
Similarly, by \eqref{jp3lp3} and \eqref{derv}
\begin{multline}\la{jl456}
  \!\!\!\!\! \Om(\tau_{j+3},\tau_{l+3})
  \! =\!-\!\!\int\!\frac{k_jk_l\Big(2i(k^2+m^2)kv(|\hat\psi_{v1}|^2\!+\!|\hat\psi_{v2}|^2)
  \!-\!((k^2+m^2)^2\!+\!(kv)^2)\hat\psi_{vc}\Big)dk}{((k^2+m^2)^2-(kv)^2)^2}=0.
\end{multline}
Finally, by (\re{Wx}),
\be\la{Wk}
   \Om(\tau_j,\tau_{l+3})
   =\int \frac{k_jk_l \Big((k^2+m^2)(|\hat\psi_{v1}|^2+|\hat\psi_{v2}|^2)+ikv\hat\psi_{vc}\Big)~dk}
   {(k^2+m^2)^2-(kv)^2} +e_j\cdot e_l.
\ee
Now (\re{Omega}) - (\re{alpha}) are proved.


\subsection*{C. Positivity of  $f_1$ and $f_2$}


Here we check the inequalities which we have used in the proof of Lemma \re{lnW}:
\begin{multline*}
1)\, f_1=\Bigl(\fr1{M-|v|k_1-\om}+\fr1{M-|v|k_1+\om} -\fr2{M-|v|k_1}\Bigr)\\
     +\Bigl(\fr1{M+|v|k_1-\om}+\fr1{M+|v|k_1+\om}-\fr2{M+|v|k_1}\Bigr)>0,
\end{multline*}
\begin{multline}\la{f2}
2)\, f_2=\Bigl(\fr1{M-|v|k_1-\om}+\fr1{M-|v|k_1+\om}-\fr2{M-|v|k_1}\Bigr)\\
     -\Bigl(\fr1{M+|v|k_1-\om}+\fr1{M+|v|k_1+\om}-\fr2{M+|v|k_1}\Bigr)\ge0
\end{multline}
under the conditions $|v|<2m,\, 0<|\om|\le \mu=m^2-v^2/4.$
First, let us note that the expressions in each bracets is positive,
since
$$\fr 1{b-a}+\fr 1{b+a}-\fr 2b=\fr{2a^2}{(b+a)(b-a)b}>0$$
if $b-a,\,b+a\ge 0,\;b>0$
and it immediately implies that $f_1>0$.
Next, the first summand in LHS of \eqref{f2} obviously is not less than the second summand
since $|v|k_1\ge0$.
Therefore $f_2\ge0$ and $f_2<f_1$.


\end{document}